\documentclass[11pt]{article}
\usepackage[authoryear]{natbib}
\usepackage[margin=25mm]{geometry}

\usepackage{graphicx}
\usepackage{subcaption}
\usepackage{float}
\usepackage{todonotes}
\usepackage{pdflscape}

\usepackage{amssymb}
\usepackage{enumitem}
\usepackage{microtype}
\usepackage{adjustbox}
\usepackage[nohyperlinks]{acronym}
\usepackage{xcolor}
\setlist{nosep}
\usepackage{hyperref}
\usepackage{amsmath}
\usepackage{mathtools}
\usepackage{tikz}
\usepackage{tikz-network}
\usepackage{textcomp}
\usepackage{multirow}
\usepackage{pifont}
\usepackage[flushleft]{threeparttable}
\newcommand{\cmark}{\ding{51}}
\newcommand{\xmark}{-}
\newcommand\Tstrut{\rule{0pt}{2.6ex}}       
\newcommand\Bstrut{\rule[-0.9ex]{0pt}{0pt}}

\definecolor{ashgrey}{rgb}{0.7, 0.75, 0.71}
\definecolor{gray(x11gray)}{rgb}{0.75, 0.75, 0.75}

    \newcommand{\stationSet}{\mathcal{S}}
    \newcommand{\transferStationSet}{\stationSet^T}

    \newcommand{\budget}[1]{b^{#1}} 
    \newcommand{\capacity}{q}
    \newcommand{\demand}[3]{d_{#1#2}^{#3}} 
    
    \newcommand{\symmetryMultiplier}{\sigma}
    \newcommand{\stationPar}{s}
    \newcommand{\firstStation}[1]{\stationPar^{#1}_{1}} 
    \newcommand{\lastStation}[1]{\stationPar^{#1}_{n_\linePar}} 
    \newcommand{\terminalStationSet}{\stationSet^{E}}
    \newcommand{\linePar}{l}
    \newcommand{\linepool}{\mathcal{L}}
    \newcommand{\stationsOnLine}[1]{\stationSet^{#1}} 
    \newcommand{\frequencySet}[1]{\mathcal{F}^{#1}} 
    \newcommand{\maxFreq}[1]{f^{#1}_{max}} 
    \newcommand{\capacityPerLine}[1]{\capacity^{#1}} 
    \newcommand{\frequencyPar}{i}
    \newcommand{\linesLeavingStation}[1]{\linepool_{#1}^{-}} 
    \newcommand{\linesEnteringStation}[1]{\linepool_{#1}^{+}} 
    \newcommand{\costOfLine}[1]{k^{#1}} 
    \newcommand{\nodeSet}{\mathcal{\text{V}}}

    \newcommand{\node}{\nu}
    \newcommand{\nodeArr}[2]{\node_{#1,#2}^{\text{Arr}}} 
    \newcommand{\nodeDep}[2]{\node_{#1,#2}^{\text{Dep}}} 
    \newcommand{\nodeIn}[1]{\node_{#1}^{\text{In}}} 
    \newcommand{\nodeOut}[1]{\node_{#1}^{\text{Out}}} 
    \newcommand{\nodeCh}[1]{\node_{#1}^{\text{Change}}} 
    \newcommand{\nodeDemand}[3]{\delta^{#2#3}_{#1}}
    \newcommand{\arcSet}[1]{\mathcal{A}_{#1}} 
    \newcommand{\arcSetDrive}[1]{\arcSet{#1}^{\text{Drive}}}
    \newcommand{\arcSetSkip}[1]{\arcSet{#1}^{\text{Skip}}}
    \newcommand{\arcSetStop}[1]{\arcSet{#1}^{\text{Stop}}}
    \newcommand{\arcSetIn}[1]{\arcSet{#1}^{\text{In}}}
    \newcommand{\arcSetOut}[1]{\arcSet{#1}^{\text{Out}}}
    \newcommand{\arcSetChangeIn}[1]{\arcSet{#1}^{\text{In-ch}}}
    \newcommand{\arcSetChangeOut}[1]{\arcSet{#1}^{\text{Out-ch}}}
    \newcommand{\arcPar}{a}
    \newcommand{\arcCost}[1]{t^{#1}}
    \newcommand{\stopVar}[3]{x^{#1,#3}_{#2}} 
    \newcommand{\frequencyVar}[3]{f^{#1,#3}_{#2}} 
    \newcommand{\flowVar}[3]{y^{#1,#3}_{#2}}

    \newcommand{\changeFreqVar}[2]{\widetilde{f}^{#1,#2}} 
    \newcommand{\changeStopVar}[3]{\widetilde{x}^{#1,#2}_{#3}}

    \newcommand{\slackVar}{e}
    \newcommand{\periodSet}{\mathcal{P}}
    \newcommand{\periodPar}[1]{p^{#1}} 
    \newcommand{\periodLength}[1]{h^{#1}} 
    \newcommand{\mop}{P1} 
    \newcommand{\mo}{P2} 
    \newcommand{\ap}{P3} 
\newcommand{\eps}{\(\varepsilon\)}

\allowdisplaybreaks

\begin{document}
\captionsetup{labelsep=space}
{\centering
{\fontsize{17.28}{16}
\textbf{Multi-period line planning for varying railway passenger demand with asymmetric lines}}}

\vspace{16pt}
\noindent R.~J.~H. {van der Knaap}\textsuperscript{a}, N. {van Oort}\textsuperscript{a}, M. {de Bruyn}\textsuperscript{b} and R.~M.~P. Goverde\textsuperscript{a}

\vspace{8pt}
\noindent\textsuperscript{a}\textit{Department of Transport and Planning, Delft University of Technology, P.O. Box 5048, 2600 GA Delft, The Netherlands}

\noindent\textsuperscript{b}\textit{Netherlands Railways, P.O. Box 2025, 3500 HA Utrecht, The Netherlands}

\vspace{16pt}
\noindent\textbf{Abstract}: A line plan is an important aspect of the quality of the service provided to railway passengers. 
Although it is well-known that railway demand is varying throughout the day in volume and structure, the line plan is often still fixed throughout the day. 
To better match this varying railway demand, we propose a mixed-integer linear programming model for multi-period line planning. 
This model for railway networks incorporates selection of routes, stopping patterns, frequencies, transfers, and the possibility of asymmetric lines to deal with spatially unbalanced demand. 
The \eps-constraint method is used to determine Pareto optimal solutions.
The proposed model and solution method are tested on a case study of part of the Dutch railway network.
The results show that allowing for changes to the line plan during the day can reduce the total generalised journey time by up to 4.26\%, especially when asymmetric lines are used.

\vspace{8pt}

\section{Introduction} \label{sec:introduction}

A railway line plan, consisting of routes, stops, and frequencies, determines a large part of the quality of the railway service for passengers.
For example, the line plan determines if passengers can have a direct trip or if they need one or multiple transfers. 
Furthermore, a line's stopping pattern influences the passengers' in-vehicle times and hence their total travel time. 
Due to this impact on passengers and the high cost of public transport networks, the line planning problem gets a lot of attention in the literature (see e.g., \cite{duran-micco2022SurveyTransitNetwork}).
Moreover, it is well-known that the passenger demand is not fixed, but instead changes throughout the day and throughout the week (see e.g., \cite{vanderknaap2024ClusteringRailwayPassenger}). 
Hence, if one single line plan is operated throughout the day, the attractiveness of this plan for the passengers can also change due to these differences in demand.
\cite{schobel2012LinePlanningPublic} points out that it is worthy to investigate whether demand in peak and off-peak hours are best served by different line plans or by the same line plan.
Moreover, \cite{zhang2020HowOptimizeTrain} show that when demand is diverse, flexible stopping rules are better able to find good solutions for the passengers compared to fixed stopping rules.
Although in urban railway it is common to have multiple schedules during the day \citep{vanoort2011ServiceReliabilityUrban}, the line plans used by railway undertakings on conventional railways in the Netherlands and other parts of Europe are still more or less fixed throughout the day.
Therefore, this paper investigates how to create a train line plan that better matches the varying demand throughout the day. 

Demand for public transport varies throughout the year (e.g., with less demand during the summer holiday as many people are on holiday), throughout the week (with different demand on different days of the week), and throughout the day \citep{vanoort2011ServiceReliabilityUrban}.
In this paper we focus on demand variations during the day.
There are several ways in which the demand can vary throughout the day. 
First, the demand volumes change: demand is higher in the peak periods at the beginning and the end of the workday.
For instance, in the Netherlands on a regular working day in 2019, around 50\% of all trips were made inside the peak hours \citep{vanderknaap2024ClusteringRailwayPassenger}.
Since the peak hours in the Netherlands are between 6:30-9:00 and 16:00-18:30, this means that the other 50\% of the trips are made during the remaining hours of the day. 
Second, the structure of the demand (i.e., which origin-destination (OD) pairs are important) also varies throughout the day.
For instance, during peak hours many people are commuting from home to their workplace or vice versa. 
Therefore, railway stations close to a business park might be highly visited during the peak hours, but not so much during working hours or in the evening. 
On the other hand, during the day there are more people travelling for leisure purposes, with for example a city centre as destination. 
Furthermore, jobs might not be equally spread over the country. 
Compared to the average, some areas have more houses and other areas have more job locations. 
During the day, this would create two distinct transport flows: in the morning more people will travel from residential areas to the areas with jobs, while in the late afternoon the direction will be reversed as people will be travelling home.

In recent years, there has been more interest in finding line plans that better match demand that is varying throughout the day.
For example, \cite{sahin2020MultiperiodLinePlanning} are the first ones to introduce a line planning model that has multiple periods during the day to accommodate varying demand and \cite{nie2023WeeklyLinePlanning} propose the Weekly Line Planning problem that looks at the demand during an entire week. 
However, there are several gaps in the multi-period line planning literature.
First, most papers on multi-period line planning use a corridor to test the proposed models. 
However, railway systems in for example Europe are typically highly connected networks.
Since it is more difficult to find a line plan for a network than for a corridor, more research is needed into how a multi-period line plan can be created for a railway network.
Second, the existing models for multi-period line planning only incorporate a limited number of the possible service adjustments to match the line plan to the demand.
Given the temporal and spatial differences in demand that are described above, possible service adjustments that can be used include changing the stopping patterns, the frequencies, the route choice, and having asymmetric lines.
We define asymmetric lines as lines that are not operated with the same stopping pattern and/or frequencies in both directions. 
Having asymmetric lines can for example better serve a demand that is not balanced in both directions.
Furthermore, flexibility in the stopping pattern in a railway network can only be achieved when passengers are able to transfer between lines. 
Therefore, including the possibility to transfer is essential for a multi-period line planning model for a railway network.  
However, to the best of our knowledge there is no model in the existing literature that incorporates transfers, frequency choice, stop choice, and asymmetric lines for multi-period line planning.

Railway undertakings can benefit from multi-period line planning, as it has the potential of better serving the passenger demand throughout the day for the same costs.
Therefore, in this paper we will fill the aforementioned gaps.
The main contribution of this paper is a novel model for multi-period line planning in a railway network that incorporates setting routes, stopping patterns, and frequencies, and includes the possibility for passengers to transfer between lines. 
Furthermore, we analyse how asymmetric lines in a railway network can improve the model's objective of minimising the passenger's generalised journey time.
The proposed model offers exceptional flexibility in aligning stopping patterns with varying demand, due to the variables used to set the stopping pattern and the inclusion of asymmetric lines.
The \eps-constraint method \citep{mavrotas2013ImprovedVersionAugmented} is used to create several line plans with varying amounts of line plan adjustments throughout the day.
To test the proposed model, a case study based on real data of part of the Dutch railway network is used.

The remainder of this paper is organised as follows.
In Section \ref{sec:literatureReview}, a literature review about (multi-period) line planning is given in order to briefly summarise previous work in this area.
Next, Section \ref{sec:mathematicalModel} introduces a mathematical model for line planning that includes choice of route, frequency and stop pattern. 
Section \ref{sec:case_study} describes the case study and Section \ref{sec:results} provides the results of applying the proposed model to this case study.
Lastly, the conclusions are provided in Section \ref{sec:conclusion}.

\section{Literature review of line planning for varying demand}\label{sec:literatureReview}
The scheduling of public transport services is often a sequential process, consisting of infrastructure planning, line planning, timetabling, vehicle scheduling, and crew scheduling \citep{ceder1986BusNetworkDesign}. 
In this paper, the focus is on the Line Planning Problem (LPP).
The output of the LPP is a line plan, which is a set of lines and their corresponding frequencies. 
Here, we define a line as a path through the public transport network and a set of served stations.

The Line Planning Problem, also called the Transit Network Design Problem, is a well-researched problem in literature.
This resulted in a rich variety of literature review papers, for example \cite{schobel2012LinePlanningPublic} and \cite{duran-micco2022SurveyTransitNetwork}. 
\cite{schobel2012LinePlanningPublic} gives an overview of the different types of line planning models found in the literature, including commonly used objectives and constraints.
One direction for future research that is given is whether peak and off-peak demand should be served by the same or different line plans.
The literature review of \cite{duran-micco2022SurveyTransitNetwork} focuses on the different versions of the Transit Network Design Problem that have been discussed in the literature. 
They conclude that technological improvements, such as improved computer capacity and data collection techniques, have facilitated and increased research interest into this topic. 
However, although this is a well-researched problem, there still exists a significant gap between theory and practice due to the simplifying assumptions used in the literature.

Most papers on the LPP have only considered constant demand. 
In the most simple cases, the demand for transport is given per edge in the public transport network, which connects two subsequent stations \citep{schobel2012LinePlanningPublic}.
Another option is that the number of passengers per origin-destination (OD) pair is given (see e.g., \cite{borndorfer2008ModelsLinePlanning}; \cite{goossens2004BranchandCutApproachSolving}).
Furthermore, the demand considered is typically for one period only.
In long distance railways, like the high speed railway network in China, trains do not run every hour. 
Therefore, the demand during the period of one day is used to determine which lines are suitable, like in \cite{fu2015HierarchicalLinePlanning}.
On the other hand, when a cyclic railway timetable is the goal, the demand is usually given for one cycle of the timetable, for example one (peak) hour. 
Papers employing this strategy include \cite{nachtigall2008SimultaneousNetworkLine} and \cite{bull2019OptimisingTravelTime}.
However, using a line plan during an entire day that is only based on one (peak) period might not be optimal for the passengers.
For example, suppose there is a station in a business park, so close to many offices but not to other facilities (e.g., shops or museums).
This station will have a lot of passengers arriving or departing during the peak hours, but not many people will want to arrive at or depart from this station during the off-peak period. 
Therefore, if a line is based on peak hour demand, and hence stops at this station, then passengers using this line in the off-peak period will have longer travel times due to the (unwanted) stop at this station. 
\cite{zhang2020HowOptimizeTrain} investigate if changing the stopping pattern of a line can be beneficial for the passengers.
They look at the stopping patterns of the trains, and find that having flexible stopping rules, which are not solely based on the type of station, is better for the passengers when demand is diverse. 
Under diverse demand, the flexible stopping rules generated a line plan with a lower average number of stops per line, which results in faster travel times for the passengers.

Besides line planning papers with fixed demand, there are also several papers that consider time-dependent demand during line planning.
These papers can be divided into two groups, depending on whether they allow the lines to change over time or not. 
The first group of papers recognises that since the demand is fluctuating throughout the day or year, a line plan should be evaluated on more than one type of demand (e.g., \cite{amiripour2014DesigningLargescaleBus}, \cite{cyril2020DemandBasedModelLine}, \cite{duran-micco2022DesigningBusLine}).
However, to create a schedule that is convenient and clear for the passengers, they want to operate the same lines with a fixed stopping pattern under each type of demand.
Meanwhile, the frequencies of the selected lines can be adjusted to better match the varying demand. 
\cite{amiripour2014DesigningLargescaleBus} consider seasonal demand variations and aim to create a bus schedule that is robust to these changes.
On the other hand, \cite{cyril2020DemandBasedModelLine} and \cite{duran-micco2022DesigningBusLine} consider demand variations during the day.  
\cite{cyril2020DemandBasedModelLine} determine for each line in the bus network a suitable frequency for both the peak and off-peak hours, while \cite{duran-micco2022DesigningBusLine} also consider the difference between morning and afternoon peak demand in their evaluation.
\cite{kaspi2013ServiceorientedLinePlanning} propose an optimisation model that creates a line plan and cyclic timetable with the aim of minimising the operational cost and the total passenger travel time. 
Time-dependent demand that changes throughout the day is used to evaluate the timetable, but as the timetable is cyclic no service-adjustments are made throughout the day.

The second group of papers makes changes to both the lines and the frequencies, instead of only changing the frequencies to match the demand. 
In these papers, the considered time interval is divided into different periods with different demand and the aim is to find a suitable set of lines and corresponding frequencies for each of these periods.
The first paper on multi-period line planning is written by \cite{sahin2020MultiperiodLinePlanning}.
They aim to find a minimal cost line plan in which the transfer of vehicles from period to period is considered. 
Continuing on this work, \cite{zhao2022FluctuatingDemandOrientedOptimization} propose a model that considers time-dependent demand, determining line frequencies, transport capacity, stop-pattern and the transfers of rolling stock in a railway corridor. 
A round heuristic algorithm is used to solve a case study based on intercity railway in China.
\cite{nie2023WeeklyLinePlanning} introduce the Weekly Line Planning (WLP) problem.
In this problem, a line plan is determined for each time block of each weekday in order to maximise the total matching utility between supply and demand and to minimise the total train operational cost.
A custom genetic algorithm is provided in order to solve two large-size case studies on the high-speed railway in China.
\cite{schiewe2023LinePlanningDifferent} highlight that unrelated line plans across different periods can lead to passenger confusion and operational challenges. 
To address this, they propose incorporating an upper bound on the dissimilarity between line plans operated within a single day in the multi-period line planning problem. 
Furthermore, the study introduces three distinct measures for assessing the (dis)similarity of these line plans.
\cite{avila-ordonez2022DesignFlexibleBus} also consider keeping changes to the line plan as small as possible, when creating a bus line plan for a major event. 
Their aim is to improve the total travel time, taking into account increased demand at the event location and congestion in that area, while changing the lines as little as possible to keep it attractive for the regular passengers.
A heuristic based on Genetic Algorithms is presented to solve this flexible bus line planning problem.

There are also several papers that combine line planning and timetabling decisions under time-dependent demand.
One group of papers only considers decisions about frequencies and/or departure times at the origin.
While these can differ throughout the considered time period, the selected or given lines are kept the same. 
To serve asymmetric time-dependent demand on a metro line, \cite{mo2021ExactMethodIntegrated} consider having different frequencies, departure times, rolling stock compositions and train speed profiles for both directions on the line.
The objective of the proposed mixed-integer non-linear programming model is to optimise both the passenger waiting and travel times and the operating costs.
\cite{li2019DemandorientedTrainServices} and \cite{yuan2022IntegratedOptimizationTrain} both want to improve passenger waiting time on a metro line by including the selection of a short-turning route in the model. 
By also including the rolling stock circulation, \cite{yuan2022IntegratedOptimizationTrain} show that using short-turning can significantly reduce the passenger waiting time compared to the timetable used in practice, while keeping the number of rolling stock the same.
Another group of papers combines timetabling with stop planning decisions.
\cite{qi2021IntegerLinearProgramming} combine stop planning and timetabling to minimise the total travel time of the trains, while serving all the demand in their preferred time slot with a direct connection. 
The resulting model is tested on the Wuhan-Guangzhou high-speed railway line in China and the model is solved using CPLEX.
\cite{zhou2023JointOptimizationTimeDependent} developed a non-linear optimisation model for finding a time-dependent line plan, timetable, and differentiated ticket prices. 
The objective of the model is to both maximise the railway undertaking's revenue and minimise the total cost of the passengers.
The model is solved using a simulated annealing algorithm. 
Moreover, in the context of congested systems, altering the stop patterns can help with reducing the travel time for passengers and the costs for the railway undertaking \citep{dong2020IntegratedOptimizationTrain,qu2023OptimizingIntegratedAllStop}, improving the fairness among waiting passengers \citep{zhao2021IntegratedApproachTrain}, and improving the safety at stations \citep{shi2023SafetyorientedTrainTimetabling}.
Contrary to the papers previously mentioned in this paragraph, which all create an acyclic timetable, \cite{zhou2023LinePlanningApproach} aim to serve time-dependent demand by using a mixture of periodic and aperiodic lines. 
Periodic lines are operated during every hour and have the same route and stopping pattern throughout the day, while aperiodic lines are added to this cyclic schedule to serve specific time-dependent demand.
The mixed-integer non-linear programming model is solved using a simulated annealing algorithm.

Reasons why the majority of the line planning papers have not considered changing demand include the unavailability of demand data and the complexity of the problem. 
Since the introduction of the smart card fare collection systems, the amount of detailed demand data has increased substantially. 
This data can for example be used to create more detailed and time-dependent OD matrices, which can be used to make service adjustments based on demand \citep{pelletier2011SmartCardData}. 
Furthermore, the LPP is a problem that is difficult to solve in general. 
In fact, most versions of the LPP are proven to be NP hard, including finding a feasible line plan \citep{bussieck1998OptimalLinesPublic}, finding a line plan with minimal costs \citep{claessens1998CostOptimalAllocation}, and finding a line plan with minimal passenger travel time \citep{schobel2005LinePlanningMinimal}. 
Therefore, several techniques are used to solve the LPP, such as making assumptions to shrink the problem size \citep{guan2006SimultaneousOptimizationTransit}, using column generation \citep{borndorfer2008ModelsLinePlanning,nachtigall2008SimultaneousNetworkLine}, and adding constraints to improve the LP relaxation \citep{borndorfer2013ConfigurationModelLine}. 
However, when the instance size gets too big, metaheuristics are needed to solve the problem. 
According to \cite{duran-micco2022SurveyTransitNetwork}, metaheuristics that are used to solve the LPP are often variations of evolutionary algorithms.
The most commonly used method is using a genetic algorithm, but other methods include NSGA-II, simulated annealing, swarm optimisation techniques, bee colony optimisation, tabu search, wolf pack search, beam search, and hyper heuristics \citep{duran-micco2022SurveyTransitNetwork}.  

\begin{table}
  \fontsize{9pt}{9pt}\selectfont
  \centering
  \begin{threeparttable}
      \centering
      \caption{Overview of the relevant literature considering line planning decisions under time-dependent demand.}
      \label{tab:literature-gap}
      \begin{tabular}{lccccccccc}
          \hline
                                                               &            &               & \multicolumn{4}{l}{\textbf{Line planning}}           &                     & \textbf{} &  \textbf{Multi-}                      \\
          \textbf{}                                            & \textbf{}  &               & \multicolumn{4}{l}{\textbf{choices included}}        & \textbf{}  & \textbf{Asym}      & \textbf{period} \\ \cline{4-7}
          \textbf{Paper}                                       & \textbf{N/C} & \textbf{Mode} & \textbf{LS} & \textbf{RS} & \textbf{SP} & \textbf{F} & \textbf{T} & \textbf{lines} & \textbf{schedule}     \\ \hline \Tstrut
          \cite{li2019DemandorientedTrainServices}             & C          & M             & \xmark      & \cmark      & \xmark      & \xmark     & \xmark              & \cmark*             & \xmark                \\
          \cite{dong2020IntegratedOptimizationTrain}           & C          & R             & \xmark      & \xmark      & \cmark      & \xmark     & \xmark              & \cmark*             & \xmark                \\
          \cite{mo2021ExactMethodIntegrated}                   & C          & M             & \xmark      & \xmark      & \xmark      & \cmark     & \xmark              & \cmark              & \xmark                \\
          \cite{qi2021IntegerLinearProgramming}                & C          & R             & \xmark      & \xmark      & \cmark      & \xmark     & \xmark              & \cmark              & \xmark                \\
          \cite{zhao2021IntegratedApproachTrain}               & C          & URT           & \xmark      & \xmark      & \cmark      & \xmark     & \xmark              & \cmark              & \cmark                \\
          \cite{nie2023WeeklyLinePlanning}                     & C          & HSR           & \cmark      & \xmark      & \xmark      & \cmark     & \cmark              & \xmark              & \cmark                \\
          \cite{yuan2022IntegratedOptimizationTrain}           & C          & M             & \xmark      & \cmark      & \xmark      & \xmark     & \xmark              & \cmark              & \xmark                \\
          \cite{zhao2022FluctuatingDemandOrientedOptimization} & C          & R             & \xmark      & \xmark      & \cmark      & \cmark     & \xmark              & \cmark              & \cmark                \\
          \cite{qu2023OptimizingIntegratedAllStop}             & C          & M             & \xmark      & \xmark      & \cmark      & \xmark     & \cmark              & \cmark*             & \xmark                \\
          \cite{shi2023SafetyorientedTrainTimetabling}         & C          & M             & \xmark      & \xmark      & \cmark      & \xmark     & \xmark              & \cmark*             & \xmark                \\
          \cite{zhou2023LinePlanningApproach}                  & C          & HSR           & \xmark      & \cmark      & \cmark      & \xmark     & \cmark              & \cmark*             & \xmark                \\
          \cite{kaspi2013ServiceorientedLinePlanning}          & N          & R             & \cmark      & \xmark      & \xmark      & \cmark     & \cmark              & \xmark              & \xmark                \\
          \cite{amiripour2014DesigningLargescaleBus}           & N          & B             & \cmark      & \xmark      & \xmark      & \cmark     & \cmark              & \xmark              & \cmark                \\
          \cite{cyril2020DemandBasedModelLine}                 & N          & B             & \xmark      & \xmark      & \xmark      & \cmark     & \xmark              & \xmark              & \cmark                \\
          \cite{sahin2020MultiperiodLinePlanning}              & N          & B, M           & \xmark      & \xmark      & \xmark      & \cmark     & \xmark              & \xmark              & \cmark                \\
          \cite{avila-ordonez2022DesignFlexibleBus}            & N          & B             & \cmark      & \xmark      & \xmark      & \xmark     & \cmark              & \xmark              & \xmark                \\
          \cite{duran-micco2022DesigningBusLine}               & N          & B             & \cmark      & \xmark      & \xmark      & \cmark     & \cmark              & \xmark              & \cmark                \\
          \cite{schiewe2023LinePlanningDifferent}              & N          & B, M           & \cmark      & \xmark      & \xmark      & \cmark     & \xmark              & \xmark              & \cmark                \\
          \cite{zhou2023JointOptimizationTimeDependent}        & N          & HSR           & \xmark      & \xmark      & \cmark      & \xmark     & \cmark              & \cmark*             & \xmark                \\
          This paper                                           & N          & R             & \xmark      & \cmark      & \cmark      & \cmark     & \cmark              & \cmark              & \cmark               \Bstrut \\ \hline
      \end{tabular}
      \fontsize{8pt}{8pt}\selectfont
      \begin{tablenotes}
          \item N/C: network (N) or corridor (C), Mode: bus (B), high-speed railway (HSR), metro (M), railway (R), urban rail transit (URT), Line plan choices: Line Selection (LS), Route Selection (RS), Stop Pattern (SP), Frequency (F).  T: transfers considered, Asym lines: asymmetric lines considered, * only one direction considered.
      \end{tablenotes}
  \end{threeparttable}
\end{table}

Table \ref{tab:literature-gap} provides an overview of the key literature that incorporates line planning decisions under time-dependent demand and shows how our paper relates to the literature.
The first three columns denote the paper, the type of network analysed, distinguishing between a corridor and a full network, and the mode considered.
The next four columns denote which line planning decisions are considered in the paper. 
``Line Selection'' (LS) refers to the selection of a set of lines with pre-defined stopping pattern from a line pool.
``Route Selection'' (RS) involves selecting a path through the network without predefined stopping patterns. 
``Stop Pattern'' (SP) indicates that the model aims to determine the stopping patterns on routes.
Lastly, ``Frequency'' (F) refers to determining the line frequencies.
The last three columns indicate: whether passenger transfers are considered, if asymmetric lines with different stopping pattern and/or frequencies in opposite directions are allowed, and whether the aim is to determine a multi-period schedule.

Although there has been a lot of research already into the LPP in general and some research into multi-period line planning problem specifically, there still exist several important research gaps. 
Firstly, most studies that consider time-dependent demand focus on bus lines or train lines in railway corridors.
While adjusting frequencies and departure times at the origin station to fluctuations in demand throughout the day is common practice in bus services (see e.g., \cite{vanoort2012ImpactSchedulingService}), this is not common in railway services.
Railway operations, in contrast to buses, face strict infrastructure constraints that complicate frequency and departure time adjustments.
Changing the departure time at the origin affects the departure time at all stations throughout the line. 
As trains could interact with other trains at each station on their route, even a small change in their departure time can create multiple conflicts due to trains claiming the same infrastructure at the same time. 
Additionally, while the railway network in China can be well represented by corridors, many European networks are highly connected networks where transferring between lines is common. 
Therefore, more research is needed into creating multi-period line plans for railway networks, where the possibility of transferring between lines is included.
Moreover, existing network-related research often only considers symmetric lines, which are either fixed or chosen from a predefined line pool. 
However, these requirements reduce the flexibility to deal with asymmetric demand.
Line pools are typically based on numerous assumptions; for instance, the Dutch railway undertaking NS categorises lines into two types: Sprinter lines that stop at all stations, and Intercity lines that stop only at the large stations. 
Additionally, some medium-large stations are served by a subset of the Intercity lines. 
However, this categorisation restricts adjustments for stopping frequency at different times, limiting the ability to accommodate dynamic demand variations effectively.
To provide more options to deal with time-varying and asymmetric demand, more flexibility in determining the stopping pattern and allowing for asymmetric lines should be considered.

In response to these gaps, this paper proposes a model for the multi-period line planning problem in railway networks.
The model's input includes a fixed number of periods that each have their own length and demand.
To deal with the temporal and spatial differences between the demands of different periods, the model considers selecting routes in the network, setting the stopping patterns and the frequencies, and having asymmetric lines. 
Since the model addresses railway networks, the possibility of transferring between lines will also be taken into account.
As shown in the last row of Table \ref{tab:literature-gap}, the proposed model fills the identified gaps in the literature, by incorporating a large set of line planning decisions in the model for multi-period line planning. 
Furthermore, we are the first to analyse the effect of incorporating asymmetric lines throughout a (railway) network.

\section{Mathematical model for multi-period line planning}\label{sec:mathematicalModel}

In this section, the mathematical model for multi-period line planning is introduced. 
The aim of the model is to provide a line plan which minimises the passengers' total Generalised Journey Time (GJT) throughout the day.
In this paper, the GJT consists of the in-vehicle time, and the waiting time before boarding a train (for both the first train and when transferring).
In Section \ref{subsec:MM_lines} we introduce the notation and variables related to train lines. 
The parameters and variables used in this model are listed in Tables \ref{tab:parametersMM_multiperiod} and \ref{tab:variables_multiperiod}, respectively. 
Next, Section \ref{subsec:MM_passengers} introduces the change-and-go network that is used to model the passenger paths and calculate the GJT.
In Section \ref{subsec:MM_modelMultiPeriod} the mathematical formulation of the multi-period line planning model is given.

\begin{table}[t]
    \centering
    \small
    \caption{Parameters used in the multi-period mathematical model.}
    \label{tab:parametersMM_multiperiod}
    \begin{tabular*}{\textwidth}[pos]{ll}
        \hline
        \multicolumn{2}{l}{\Tstrut Parameters \Bstrut} \\
        \hline \Tstrut
        \(\stationSet\)             & Set of stations \\
        \(\transferStationSet\)     & Set of stations at which a transfer is allowed, \(\transferStationSet \subseteq \stationSet\) \\
        \(\terminalStationSet\)     & Set of terminal stations. \\
        \(\linepool\)               & Set of candidate lines \\
        \(\linesLeavingStation{\stationPar}\) & Set of candidate lines that have station \(\stationPar \in \terminalStationSet\) as first station, \(\linesLeavingStation{\stationPar} \subseteq \linepool\). \\
        \(\linesEnteringStation{\stationPar}\) & Set of candidate lines that have station \(\stationPar \in \terminalStationSet\) as last station, \(\linesEnteringStation{\stationPar} \subseteq \linepool\). \\
        \(\stationPar^\linePar_i\)  & The \(i\textsuperscript{th}\) station on candidate line \(\linePar \in \linepool\). \\
        \(\firstStation{\linePar}\) & The first station on candidate line \(\linePar \in \linepool\). \\
        \(\lastStation{\linePar}\)  & The last station on candidate line \(\linePar \in \linepool\). The number of stations on a candidate  \\
                                    & line \(n_\linePar\) can be different for each line. \\ 
        \(\stationsOnLine{\linePar}\)& Sequence of adjacent stations \((\firstStation{\linePar}, \dots, \lastStation{\linePar})\) defining the route of candidate line \(\linePar \in \linepool\).\\ 
        \(\frequencySet{\linePar}\) & Set of allowed frequencies for candidate line \(\linePar \in \linepool\). The frequency is the number of  \\
                                    & trains per hour operated on the line, so all frequencies are positive integers. \\
        \(\maxFreq{\linePar}\)      & The largest frequency in \(\frequencySet{\linePar}\). \\
        \(\periodSet\)              & Set of periods for which we want to determine a line plan. Let \(\periodPar{j}\) denote the \\
                                    & \(j\textsuperscript{th}\) period in \(\periodSet\). \\
        \(\periodLength{\periodPar{}}\) & Length of period \(\periodPar{} \in \periodSet\) in hours. \\
        \(\budget{\periodPar{}}\)   & The maximum budget for a line plan in period \(\periodPar{}\), given as the total number of  \\
                                    &  kilometres that trains can drive in one hour.\\
        \(\capacityPerLine{\linePar}\)& The capacity of a train on candidate line \(\linePar\). \\
        \(\demand{\stationPar}{u}{\periodPar{}}\) & Passenger demand from station \(\stationPar\) to station \(u\) (\(\stationPar,u \in \stationSet\)) in period \(\periodPar{} \in \periodSet\).\\
        \(\nodeDemand{\node}{\stationPar}{\periodPar{}}\) & Demand at vertex \(\node \in \nodeSet{}\) for passengers who originate from station \(\stationPar \in \stationSet\) in \\
                                    & period \(\periodPar{} \in \periodSet\). \\
        \(\arcCost{\arcPar}\)       & Time needed to travel on arc \(\arcPar \in \arcSet{}\) in the change-and-go network \\
        \(\costOfLine{\linePar}\)   & Length of candidate line \(\linePar\) in kilometres. \\
        \(\symmetryMultiplier\)     & Symmetry parameter, takes value 2 if symmetric lines are considered, 1 otherwise. \Bstrut \\ 
        \hline
    \end{tabular*}
\end{table}
\begin{table}[ht]
    \centering
    \small
    \caption{Variables used in the multi-period mathematical model.}
    \label{tab:variables_multiperiod}
    \begin{tabular}{ll}
        \hline
        \multicolumn{2}{l}{\Tstrut Variables} \Bstrut \\
        \hline 
        \Tstrut \(\stopVar{\linePar}{\stationPar}{\periodPar{}}\)& Binary variable which is equal to 1 if line \(\linePar\) stops at station $\stationPar$ in period \(\periodPar{}\), 0 \\
                        &   otherwise. This variable is defined for each candidate line \(\linePar \in \linepool\), for each station    \\
                        & \(\stationPar \in \stationsOnLine{\linePar} \setminus \{\firstStation{\linePar}, \lastStation{\linePar}\}\), and for each period \(\periodPar{} \in \periodSet\).\\
        \(\frequencyVar{\linePar}{\frequencyPar}{\periodPar{}}\) & Binary variable which is equal to 1 if line $\linePar$ is operated with frequency $\frequencyPar$ in period \\
                        &  $\periodPar{}$, 0 otherwise. This variable is defined for each candidate line \(\linePar \in \linepool\), each frequency \\
                        & \(\frequencyPar \in \frequencySet{\linePar}\), and each period \(\periodPar{} \in \periodSet\).\\
        \(\flowVar{\arcPar}{\stationPar}{\periodPar{}}\) & Real variable denoting the number of passengers that originate from station \(\stationPar\) and \\
                        &  use arc \(\arcPar\) to travel to their destination in period \(\periodPar{}\). This variable is defined for each \\
                        &  origin station \(\stationPar \in \stationSet\), each arc \(\arcPar \in \arcSet{}\) in the change-and-go network, and each period \\
                        & \(\periodPar{} \in \periodSet\).\\
        \(\changeFreqVar{\linePar}{\periodPar{}}\) & Real variable denoting whether the frequency of line \(\linePar\) has changed in period \(\periodPar{}\)\\ 
                        &  compared to the previous period. This variable is defined for each candidate \\
                        & line \(\linePar \in \linepool\), and each period except the first one: \(\periodPar{} \in \periodSet \setminus \{\periodPar{1}\}\). \\
        \(\changeStopVar{\linePar}{\periodPar{}}{\stationPar}\) & Real variable denoting whether the stop plan of line \(\linePar\) at station \(\stationPar\) has changed in \\
                        & period \(\periodPar{}\) compared to the previous period. This variable is defined for each candidate  \\
                        & line \(\linePar \in \linepool\), for each station \(\stationPar \in \stationsOnLine{\linePar} \setminus \{\firstStation{\linePar}, \lastStation{\linePar}\}\), and each period except the first one:  \\
                        & \(\periodPar{} \in \periodSet \setminus \{\periodPar{1}\}\). \Bstrut \\
        \hline
    \end{tabular}
\end{table}

\subsection{Modelling train lines}\label{subsec:MM_lines}
To minimise the passengers' total GJT when the demand varies throughout the day, the model should be able to create different lines for different periods.
A train line is defined by its route through the network, its stopping pattern, and its frequency (i.e., how many times per hour a line is operated).
Contrary to most literature, our proposed model does not select lines with fixed route and stopping pattern from a line pool. 
Instead, similarly to \cite{zhou2023JointOptimizationTimeDependent} each candidate line has a fixed route, but its stopping pattern and frequency are determined by the model using decision variables.
As part of the input, we need to provide a set of routes and per route how many candidate lines should be included.
A route is a sequence of adjacent stations in the railway network, that determines the path through the railway network.
We say that two stations are adjacent if they are connected by a track.
The first and last station of the sequence are the two terminal stations.
Usually there is only a limited number of stations that can function as a terminal station.
To qualify as a terminal station, a station should for example have infrastructure for turning the trains and/or capacity for parking rolling stock that is temporarily not used in service.
Hence, the selection of which stations are terminal stations greatly influences which routes can be included in the route set.
Furthermore, as it could be desirable to have multiple lines with a different stopping pattern on the same route, we allow for more than one candidate line per route.
The set of candidate lines is denoted by \(\linepool\).

The stopping pattern of a candidate line is set by the stopping variables.
Each station \(\stationPar\) on candidate line \(\linePar\) gets a binary stopping variable \(\stopVar{\linePar}{\stationPar}{\periodPar{}}\), which takes value one if the line stops at station \(\stationPar\) in period \(\periodPar{}\) and zero otherwise.
We assume that a line has to stop at the terminal stations, so no stopping variable is created for these stations.

Besides stopping variables, candidate lines also have frequency variables.
Each candidate line \(\linePar\) has a set of frequencies \(\frequencySet{\linePar}\), denoting how many trains per hour are allowed to be operated on the line.
For each of these frequencies, a binary decision variable \(\frequencyVar{\linePar}{\frequencyPar}{\periodPar{}}\) is created that is equal to one if the line is operated with frequency \(\frequencyPar \in \frequencySet{\linePar}\) in period \(\periodPar{}\) and zero otherwise.
If a candidate line is not selected in the final line plan, all its frequency variables will take value zero. 

Rolling stock allocation is not taken into account in this model.
Instead, for each route a fixed train capacity is given in the input.
This capacity could for example be the capacity of the largest train that can be operated on that route.

\subsection{Modelling passengers and travel time}\label{subsec:MM_passengers}
The objective of our mathematical model is to minimise the total Generalised Journey Time of the passengers. 
The GJT was first introduced by \cite{tyler1973GravityElasticityModels} to evaluate the timetable related attractiveness of different train services.
To calculate the in-vehicle time and waiting times for the GJT, we assign the passengers to arcs in the change-and-go (C\&G) network.
The C\&G network was first introduced by \cite{schobel2005LinePlanningMinimal} to model passenger transfers in a public transport network.
In this paper we use two versions of the C\&G network: one with asymmetric lines  (Section \ref{subsec:model_asym_lines}) and one with symmetric lines (Section \ref{subsubsec:model_sym_lines}).

\subsubsection{Modelling passengers with asymmetric lines}\label{subsec:model_asym_lines}
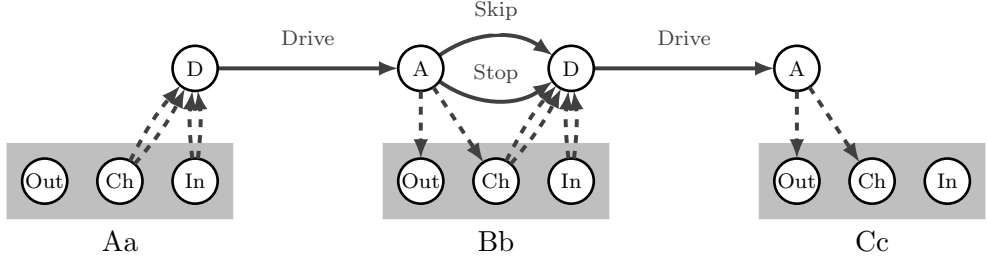
\begin{figure}[t]
    \centering
    \begin{tikzpicture}
        \filldraw[gray(x11gray)] (5.5,2) rectangle (8.5,3);
        \filldraw[gray(x11gray)] (0.5,2) rectangle (3.5,3);
        \filldraw[gray(x11gray)] (10.5,2) rectangle (13.5,3);
        \node at (2,1.7) {Aa};
        \node at (7,1.7) {Bb};
        \node at (12,1.7) {Cc};
        \Vertex[x=6,y=4,label=A,color=white]{A1}
        \Vertex[x=11,y=4,label=A,color=white]{A2}
        \Vertex[x=3,y=4,label=D,color=white]{D1}
        \Vertex[x=8,y=4,label=D,color=white]{D2}
        \Vertex[x=8,y=2.5,label=In,color=white]{In}
        \Vertex[x=7,y=2.5,label=Ch,color=white]{Ch}
        \Vertex[x=6,y=2.5,label=Out,color=white]{Out}
        \Vertex[x=3,y=2.5,label=In,color=white]{In2}
        \Vertex[x=2,y=2.5,label=Ch,color=white]{Ch2}
        \Vertex[x=1,y=2.5,label=Out,color=white]{Out2}
        \Vertex[x=13,y=2.5,label=In,color=white]{In3}
        \Vertex[x=12,y=2.5,label=Ch,color=white]{Ch3}
        \Vertex[x=11,y=2.5,label=Out,color=white]{Out3}
        \Edge[Direct,label=Drive,position=above](D1)(A1)
        \Edge[Direct,label=Drive,position=above](D2)(A2)
        \Edge[Direct,bend=35,label=Skip,position=above](A1)(D2)
        \Edge[Direct,bend=-35,label=Stop,position=above](A1)(D2)
        \Edge[Direct,NotInBG=true,style=dashed](A1)(Out)
        \Edge[Direct,NotInBG=true,style=dashed](A1)(Ch)
        \Edge[Direct,NotInBG=true,style=dashed,bend=-7](Ch)(D2)
        \Edge[Direct,NotInBG=true,style=dashed,bend=7](Ch)(D2)
        \Edge[Direct,NotInBG=true,style=dashed,bend=-7](In)(D2)
        \Edge[Direct,NotInBG=true,style=dashed,bend=7](In)(D2)
        \Edge[Direct,NotInBG=true,style=dashed](A2)(Out3)
        \Edge[Direct,NotInBG=true,style=dashed](A2)(Ch3)
        \Edge[Direct,NotInBG=true,style=dashed,bend=-7](Ch2)(D1)
        \Edge[Direct,NotInBG=true,style=dashed,bend=7](Ch2)(D1)
        \Edge[Direct,NotInBG=true,style=dashed,bend=-7](In2)(D1)
        \Edge[Direct,NotInBG=true,style=dashed,bend=7](In2)(D1)
    \end{tikzpicture}
    \caption{Example of a change-and-go network with one line that is operated in one direction.}
    \label{fig:C&Gnetwork}        
\end{figure}
In this section, we describe the C\&G network that is used when asymmetric lines are considered.
In Figure \ref{fig:C&Gnetwork}, a small example of the C\&G network is given.
This example contains three stations (Aa, Bb, and Cc) and one (one-directional) train line from station Aa to station Cc that can be operated at two different frequencies. 
The different types of nodes and arcs are introduced in the next two paragraphs.

In the C\&G network, each train line gets a Departure (D) node at each station on the line except the destination station and an Arrival (A) node at each station except the first station on the line. 
Each station in the C\&G network is represented by three nodes: the In, Change, and Out node.
Passengers can access the public transport network at the In node and leave the network at the Out node.
The Change node can be used to change lines at that station. 

The C\&G network also contains seven types of arcs connecting the different nodes.
In Figure \ref{fig:C&Gnetwork}, the arcs related to train lines are denoted by solid lines, while the arcs used by passengers to board and alight trains are denoted by dashed lines.
There are three types of arcs related to train lines: Drive arcs, Skip arcs and Stop arcs.
Drive arcs are added from Departure nodes to the Arrival nodes at the next station and represent the train driving to the next station on the line.
The Skip and Stop arcs are added at each intermediate station on the line.
If the train line stops at the intermediate station, the Stop arc will be used, and otherwise the Skip arc will be used. 
Next, there are four types of arcs that are related to passengers: In, Out, In-change, and Out-change arcs.
From each Arrival node, there is an Out arc going to the Out node of the station for passengers that arrived at their destination and leave the network.
The Out-change arc is created from the Arrival node to the Change node, for people who would like to change lines at this station.
The In and In-change arcs are arcs from the In and Change nodes of the station, respectively, to the Departure node of the train line. 
Following the work of \cite{bull2019OptimisingTravelTime}, instead of one arc from the In and Change node to the train line, we introduce one arc per possible frequency of the train line. 
In the example of Figure \ref{fig:C&Gnetwork}, the train line can have two different frequencies and therefore there are two arcs going from each In and Change node to the Departure nodes of the train line. 
By introducing an arc for boarding the train for each frequency, we can set the cost of those arcs to a frequency dependent waiting time for boarding a line without introducing quadratic terms in our mathematical model.

The notation of the C\&G network is given in Table \ref{tab:notationCnGnetwork}.
All arcs in the C\&G network are linked to a certain candidate line: either arcs are used to describe a movement of a train on the line (Drive, Skip or Stop), or arcs are used to describe movements of passengers boarding or alighting trains on the line (In(-change) and Out(-change)).
Furthermore, we can link all arcs, except the Drive arcs, to a specific station, because a line  stops at or skips a station and passengers board or alight at a station.
It is useful to have for each type of arc (e.g., Stop arcs), except Drive arcs, a separate set of arcs connected to a line \(\linePar\) (denoted by \(\arcSetStop{\linePar}\)) and a set of arcs connected to a station \(\stationPar\) (denoted by \(\arcSetStop{\stationPar}\)).
Then, if we want to create a constraint for a Stop arc of candidate line \(\linePar\) at station \(\stationPar\), we can find the correct arc by taking the intersection of the sets: \(\arcSetStop{\linePar} \cap \arcSetStop{\stationPar}\).
\begin{table}[t]
    \centering
    \small
    \caption{Notation change-and-go network.}
    \label{tab:notationCnGnetwork}
    \begin{tabular}{ll}
        \hline \Tstrut
        \(\nodeSet{}\)              & Set of all nodes in the change-and-go network. \\
        \(\nodeArr{\linePar}{\stationPar}\) & Denotes the  arrival of line \(\linePar\) at station \(\stationPar \in \stationsOnLine{\linePar}\). This node is generated for all \(\linePar \in \linepool\)\\
                                    &  and for all stations \(\stationPar \in \stationsOnLine{\linePar}\), \(\stationPar \neq \firstStation{\linePar}\).\\
        \(\nodeDep{\linePar}{\stationPar}\) & Denotes the departure of line \(\linePar\) at station \(\stationPar \in \stationsOnLine{\linePar}\). This node is generated for all\\
                                    &  \(\linePar \in \linepool\) and for all stations \(\stationPar \in \stationsOnLine{\linePar},\stationPar \neq \lastStation{\linePar}\). \\
        \(\nodeIn{\stationPar}\)    & The source node of station \(\stationPar\), which is generated for each \(\stationPar \in \stationSet\). \\
        \(\nodeOut{\stationPar}\)   & The sink node of station \(\stationPar\), which is generated for each \(\stationPar \in \stationSet\). \\
        \(\nodeCh{\stationPar}\)    & The change node of station \(\stationPar\), which is generated for each \(\stationPar \in \transferStationSet\). \\
        \(\arcSet{}\)               & Set of all arcs in the change-and-go network \\
        \(\arcSetDrive{\linePar}\)  & Set of arcs which represent driving from one station to the next on line \(\linePar \in \linepool\). A \\
                                    & drive arc \((\nodeDep{\linePar}{\stationPar_j^{\linePar}}, \nodeArr{\linePar}{\stationPar_{j+1}^{\linePar}})\) is created for each candidate line \(\linePar \in \linepool\) and for each \\
                                    & station \(\stationPar_j^{\linePar} \in \stationsOnLine{\linePar} \setminus \lastStation{\linePar}\).\\
        \(\arcSetSkip{\linePar}\)   & Set of arcs which represent skipping a stop on line \(\linePar \in \linepool\). A skip arc \((\nodeArr{\linePar}{\stationPar}, \nodeDep{\linePar}{\stationPar})\) \\ 
                                    & is created for each candidate line \(\linePar \in \linepool\) and each station \(\stationPar \in \stationsOnLine{\linePar} \setminus \{\firstStation{\linePar}, \lastStation{\linePar}\}\).\\
        \(\arcSetStop{\linePar}\)   & Set of arcs representing a stop at a station on line  \(\linePar \in \linepool\). A stop arc \((\nodeArr{\linePar}{\stationPar}, \nodeDep{\linePar}{\stationPar})\)\\
                                    &   is created for each candidate line \(\linePar \in \linepool\) and each station \(\stationPar \in \stationsOnLine{\linePar} \setminus \{\firstStation{\linePar}, \lastStation{\linePar}\}\). \\
        \(\arcSetIn{\linePar}\)     & Set of arcs representing passengers boarding line \(\linePar \in \linepool\) from their origin station.\\
                                    & For each candidate line \(\linePar \in \linepool\), for each \(\stationPar \in \stationsOnLine{\linePar}\setminus \lastStation{\linePar}\), and for each frequency \\
                                    & \(i \in \frequencySet{\linePar}\) an entrance arc \((\nodeIn{\stationPar}, \nodeDep{\linePar}{\stationPar})\) is created.\\
        \(\arcSetOut{\linePar}\)    & Set of arcs which represent passengers leaving line \(\linePar \in \linepool\) because they arrived at \\
                                    & their destination station. For each candidate line \(\linePar \in \linepool\) and for each \(\stationPar \in \stationsOnLine{\linePar}\setminus \firstStation{\linePar}\) \\
                                    & an exit arc \((\nodeArr{\linePar}{\stationPar}, \nodeOut{\stationPar})\) is created. \\
        \(\arcSetChangeIn{\linePar}\)& Set of arcs which represent passengers entering line \(\linePar \in \linepool\) after transferring at a \\
                                    &  station. For each candidate line \(\linePar \in \linepool\), for each station  \(\stationPar \in (\stationsOnLine{\linePar}\cap \transferStationSet) \setminus \lastStation{\linePar}\), and for  \\
                                    & each frequency \(i \in \frequencySet{\linePar}\) an incoming change arc \((\nodeCh{\stationPar}, \nodeDep{\linePar}{\stationPar}; i)\) is created. \\
        \(\arcSetChangeOut{\linePar}\)& Set of arcs which represent passengers leaving line \(\linePar \in \linepool\) to  transfer at a station. \\
                                    & For each candidate line \(\linePar \in \linepool\) and for each station \(\stationPar \in (\stationsOnLine{\linePar}\cap \transferStationSet) \setminus \firstStation{\linePar}\) an outgoing \\
                                    & change arc \((\nodeArr{\linePar}{\stationPar}, \nodeCh{\stationPar})\) is created. \\
        \(\arcSet{}^{X}\)           & Set of all arcs of type \(X\), where \(X \in\) \{Drive, Skip, Stop, In, Out, Ch-in, Ch-out\}. \\
                                    & \(\arcSet{}^{X}\coloneqq \cup_{\linePar \in \linepool} \arcSet{\linePar}^{X}\). \\
        \(\arcSet{\stationPar}^{X}\)& Set of arcs of type \(X\) linked to station \(\stationPar \in \stationSet\), where \(X \in\) \{Skip, Stop, In, Out, \\
                                    & Ch-in, Ch-out\}. Note that \(\arcSet{\stationPar}^{X} \subset \arcSet{}^{X}\). \\
        \hline
    \end{tabular}
\end{table}

The GJT can be calculated using the C\&G network, by assigning a travel time related cost to each arc and determining the number of passengers that use that arc. 
We use \(\arcCost{\arcPar}\) to denote the cost of using arc \(\arcPar\) and this cost is assumed to be given. 
For each type of arc, a different cost is set. 
When \(\arcPar\) is a Drive arc \((\arcPar \in \arcSetDrive{})\), \(\arcCost{\arcPar}\) is the driving time between stations. 
For Stop arcs \((\arcPar \in \arcSetStop{})\), \(\arcCost{\arcPar}\) represents the dwell time at a station and a correction for the time loss due to braking before the stop and accelerating after the stop. 
If \(\arcPar \in \arcSetSkip{}\), \(\arcCost{\arcPar}\) is set to 0 to represent the absence of time loss by skipping the stop. 
For Out and Out-change arcs \((\arcPar \in \arcSetOut{} \cup \arcSetChangeOut{})\), \(\arcCost{\arcPar}\) is set to a fixed penalty that represents the time-loss due to braking.
Lastly, if \(\arcPar \in \arcSetIn{} \cup \arcSetChangeIn{}\), the arc cost will be a frequency-dependent penalty for the waiting time plus a fixed correction for acceleration of the train.

To track the routing of passengers through the C\&G network, passenger flow variables are used.
Variable \(\flowVar{\arcPar}{\stationPar}{\periodPar{}}\) denotes the number of passengers that originate from station \(\stationPar\) and use arc \(\arcPar\) in period \(\periodPar{}\).
This grouping of passenger flows by origin station was introduced by \cite{bull2019OptimisingTravelTime} and significantly reduces the number of variables needed to represent the passenger flows, compared to having a flow variable for each origin-destination pair. 
To ensure that all passengers are assigned a route through the network from their origin to their destination, the passenger demand is transformed to node demand.
Following the work of \cite{bull2019OptimisingTravelTime}, let the demand at vertex \(\node\) for passengers who originate from station \(\stationPar\) in period \(\periodPar{}\) be denoted by \(\nodeDemand{\node}{\stationPar}{\periodPar{}}\).
This demand represents the number of passengers entering the node minus the number of passengers leaving the node, and is defined as follows:
\begin{equation*}
    \nodeDemand{\node}{\stationPar}{\periodPar{}} = 
    \begin{cases}
        \demand{\stationPar}{t}{\periodPar{}} & \text{if vertex $\node$ is the Out node for station $t$,} \\
        - \sum_{t \in \stationSet} \demand{\stationPar}{t}{\periodPar{}} & \text{if vertex $\node$ is the In node for station } \stationPar, \\
        0 & \text{otherwise.}
    \end{cases}
\end{equation*}
If this demand is satisfied at each node in the C\&G network, all passengers travel from their origin to their destination and the passenger flow is conserved when a node is not the origin or the destination.
With the flow variables, we can determine how many passengers use each arc in the C\&G network. 
Therefore, we can calculate the total GJT by multiplying the cost of each arc with the number of passengers using that arc and summing over all arcs:
\begin{equation}
    \sum_{\periodPar{} \in \periodSet} \sum_{\arcPar \in \arcSet{}} \sum_{\stationPar \in \stationSet} \periodLength{\periodPar{}} \arcCost{\arcPar} \flowVar{\arcPar}{\stationPar}{\periodPar{}}.\label{MM-MP:obj_GJT}
\end{equation}
As each period can have a different length, the calculated GJT for each period is multiplied by the length of that period \((\periodLength{\periodPar{}})\) to get a fair comparison.

\subsubsection{Adaption for modelling symmetric lines}\label{subsubsec:model_sym_lines}
The previous section assumes that the candidate lines are defined in only one direction. 
Since each candidate line has its own route and variables for frequency and stop choice, the choices about stopping pattern and frequencies of two lines whose routes pass the same stations but in opposite directions, are independent of each other.
This gives the model more opportunities to reduce the total GJT, as it can make different choices about the stopping patterns and frequencies for opposite directions when the demand on a route is asymmetric.
However, currently in the Netherlands the line plan only contains symmetric lines. 
A symmetric line is operated in both directions with the same frequency and stopping pattern, which makes a line plan with only symmetric lines easier to understand and remember for the passengers.
Hence, we also create a version of the model that has symmetric lines, such that we can see what the benefits and drawbacks are from having only symmetric lines.

One way to ensure the creation of symmetric lines is by adding constraints stating that the value of the stop variables and frequency variables should be equal for two lines that run in opposite directions.  
However, by changing the C\&G network instead, the number of nodes and arcs in the C\&G network can be reduced as well as the number of variables and constraints needed.  
As this has a positive effect on the scalability and solvability of the model, we choose for this method.
In the remainder of this section, we describe which changes should be made to the C\&G network when we assume that lines are operated in two directions instead of one.

In Figure \ref{fig:C&Gnetwork_symmetricLines}, an example of a C\&G network with a symmetric line is shown.
The difference between this network and the C\&G network with asymmetric lines, as displayed in Figure \ref{fig:C&Gnetwork}, consists of the following changes:
\begin{itemize}
    \item An Arrival and Departure node is created for every station on the line. In the asymmetric case, no arrival node is created for the first station on the line and no departure node is created for the last station on the line.
    \item Drive arcs are added in both directions, instead of in only one direction. 
    \item The In and In-change arcs are created for all stations. In the asymmetric case, these arcs are not created for the last station on the line. 
    \item The Out and Out-change arcs are created for all stations. In the asymmetric case, these arcs are not created for the first station on the line. 
\end{itemize}
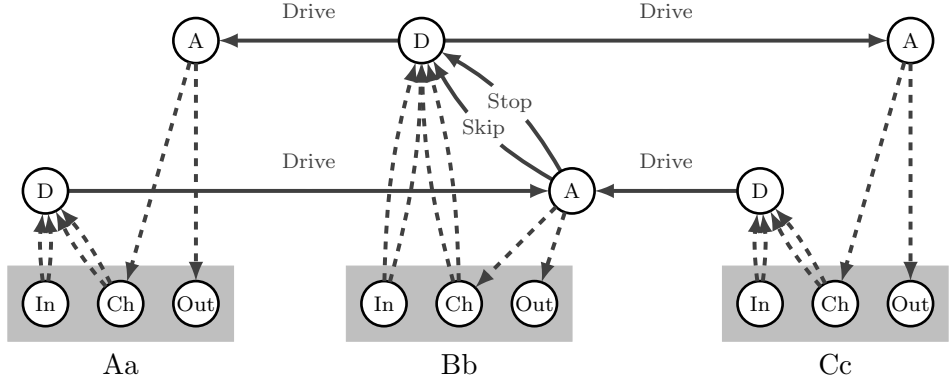
\begin{figure}[t]
    \centering
    \begin{tikzpicture}
        \filldraw[gray(x11gray)] (5,2) rectangle (8,3);
        \filldraw[gray(x11gray)] (0.5,2) rectangle (3.5,3);
        \filldraw[gray(x11gray)] (10,2) rectangle (13,3);
        \node at (2,1.7) {Aa};
        \node at (6.5,1.7) {Bb};
        \node at (11.5,1.7) {Cc};
        \Vertex[x=3,y=6,label=A,color=white]{A1}
        \Vertex[x=1,y=4,label=D,color=white]{D1}
        \Vertex[x=8,y=4,label=A,color=white]{A2}
        \Vertex[x=6,y=6,label=D,color=white]{D2}
        \Vertex[x=12.5,y=6,label=A,color=white]{A3}
        \Vertex[x=10.5,y=4,label=D,color=white]{D3}
        \Vertex[x=1,y=2.5,label=In,color=white]{In}
        \Vertex[x=2,y=2.5,label=Ch,color=white]{Ch}
        \Vertex[x=3,y=2.5,label=Out,color=white]{Out}
        \Vertex[x=5.5,y=2.5,label=In,color=white]{In2}
        \Vertex[x=6.5,y=2.5,label=Ch,color=white]{Ch2}
        \Vertex[x=7.5,y=2.5,label=Out,color=white]{Out2}
        \Vertex[x=10.5,y=2.5,label=In,color=white]{In3}
        \Vertex[x=11.5,y=2.5,label=Ch,color=white]{Ch3}
        \Vertex[x=12.5,y=2.5,label=Out,color=white]{Out3}
        \Edge[Direct,label=Drive,position=above](D1)(A2)
        \Edge[Direct,label=Drive,position=above](D2)(A3)
        \Edge[Direct,label=Drive,position=above](D3)(A2)
        \Edge[Direct,label=Drive,position=above](D2)(A1)
        \Edge[Direct,bend=15,label=Skip](A2)(D2)
        \Edge[Direct,bend=-15,label=Stop](A2)(D2)
        \Edge[Direct,NotInBG=true,style=dashed](A1)(Out)
        \Edge[Direct,NotInBG=true,style=dashed](A1)(Ch)
        \Edge[Direct,NotInBG=true,style=dashed,bend=-7](Ch)(D1)
        \Edge[Direct,NotInBG=true,style=dashed,bend=7](Ch)(D1)
        \Edge[Direct,NotInBG=true,style=dashed,bend=-7](In)(D1)
        \Edge[Direct,NotInBG=true,style=dashed,bend=7](In)(D1)
        \Edge[Direct,NotInBG=true,style=dashed](A2)(Out2)
        \Edge[Direct,NotInBG=true,style=dashed](A2)(Ch2)
        \Edge[Direct,NotInBG=true,style=dashed,bend=-7](Ch2)(D2)
        \Edge[Direct,NotInBG=true,style=dashed,bend=7](Ch2)(D2)
        \Edge[Direct,NotInBG=true,style=dashed,bend=-7](In2)(D2)
        \Edge[Direct,NotInBG=true,style=dashed,bend=7](In2)(D2)
        \Edge[Direct,NotInBG=true,style=dashed](A3)(Out3)
        \Edge[Direct,NotInBG=true,style=dashed](A3)(Ch3)
        \Edge[Direct,NotInBG=true,style=dashed,bend=-7](Ch3)(D3)
        \Edge[Direct,NotInBG=true,style=dashed,bend=7](Ch3)(D3)
        \Edge[Direct,NotInBG=true,style=dashed,bend=-7](In3)(D3)
        \Edge[Direct,NotInBG=true,style=dashed,bend=7](In3)(D3)
    \end{tikzpicture}
    \caption{Example of a change-and-go network with one symmetric line that is operated in both directions.}
    \label{fig:C&Gnetwork_symmetricLines}        
\end{figure}
By making these changes to the creation of the C\&G network, we can use all arcs except the Drive arcs for travelling in both directions on the same line. 
Hence, we get a C\&G with fewer arcs, which also means that fewer flow variables are needed as these variables are introduced for each arc.
Besides fewer flow variables, the number of frequency and stop variables needed can also be halved since these now determine the frequency and stops in two directions.
In order to account for these differences in the C\&G network in the mathematical model, we introduce symmetry parameter \(\symmetryMultiplier\).
Parameter \(\symmetryMultiplier\) takes value 2 if symmetric lines are considered and 1 otherwise.
In the next section, we describe how this parameter is used in the mathematical model.

\subsection{Mathematical formulation}\label{subsec:MM_modelMultiPeriod}
In this section, the parameters, variables, and notation described in Tables \ref{tab:parametersMM_multiperiod}, \ref{tab:variables_multiperiod}, and \ref{tab:notationCnGnetwork} are used to define the constraints and objective of the multi-period line planning model.

\subsubsection{Constraints}
The constraints of the multi-period line planning model are given in \eqref{MM-MP:flowConstraints}-\eqref{MM-MP:freqChangeSign}.

\begin{align}
    & \sum_{(u,v) \in \arcSet{}} \flowVar{(u,v)}{\stationPar}{\periodPar{}} - \sum_{(v,w) \in \arcSet{}} \flowVar{(v,w)}{\stationPar}{\periodPar{}} =  \!\!\!\!\!\!\!\! &&\nodeDemand{\node}{\stationPar}{\periodPar{}} \qquad \forall \stationPar \in \stationSet,\ \node \in \nodeSet,\ \periodPar{} \in \periodSet \label{MM-MP:flowConstraints} \\
    & \sum_{\frequencyPar \in \frequencySet{\linePar}} \frequencyVar{\frequencyPar}{\linePar}{\periodPar{}} \leq 1 \qquad &&\forall \linePar \in \linepool,\ \periodPar{} \in \periodSet \label{MM-MP:maxOneFrequency}\\
    & \sum_{\stationPar \in \stationSet} \flowVar{\arcPar}{\stationPar}{\periodPar{}} \leq \capacityPerLine{\linePar} \cdot \Bigl(\sum_{\frequencyPar \in \frequencySet{\linePar}} \frequencyPar \cdot \frequencyVar{\frequencyPar}{\linePar}{\periodPar{}}\Bigr)\quad &&\forall \linePar \in \linepool,\ \arcPar \in \arcSetDrive{\linePar},\ \periodPar{} \in \periodSet \label{MM-MP:useDriveArc} \\
    & \sum_{\stationPar \in \stationSet} \flowVar{\arcPar}{\stationPar}{\periodPar{}} \leq \maxFreq{\linePar} \cdot \stopVar{\linePar}{u}{\periodPar{}} \cdot \capacityPerLine{\linePar} \cdot \symmetryMultiplier  \quad &&\forall \linePar \in \linepool,\ u \in \stationsOnLine{\linePar} \setminus \{\firstStation{\linePar}, \lastStation{\linePar}\}, \arcPar \in \arcSetStop{\linePar} \cap \arcSetStop{u},\ \periodPar{} \in \periodSet \label{MM-MP:useStopArc}\\
    & \sum_{\stationPar \in \stationSet} \flowVar{\arcPar}{\stationPar}{\periodPar{}} \leq \maxFreq{\linePar} \cdot (1-\stopVar{\linePar}{u}{\periodPar{}}) \cdot \capacityPerLine{\linePar} \cdot \symmetryMultiplier \!\!\!\! &&\forall \linePar \in \linepool,\ u \in \stationsOnLine{\linePar} \setminus \{\firstStation{\linePar}, \lastStation{\linePar}\},\  \arcPar \in \arcSetSkip{\linePar} \cap \arcSetSkip{u},\ \periodPar{} \in \periodSet \label{MM-MP:useSkipArc}\\
    & \sum_{\stationPar \in \stationSet} \flowVar{\arcPar}{\stationPar}{\periodPar{}} \leq \maxFreq{\linePar} \cdot \stopVar{\linePar}{u}{\periodPar{}} \cdot \capacityPerLine{\linePar} \cdot \symmetryMultiplier  \quad &&\forall \linePar \in \linepool,\ u \in \stationsOnLine{\linePar} \setminus \{\firstStation{\linePar}, \lastStation{\linePar}\},\ \periodPar{} \in \periodSet,  \label{MM-MP:useEnterOrLeaveLineArcs}\vspace{-1em}\\
    & && \; \arcPar \in (\arcSetOut{\linePar} \cap \arcSetOut{u})\ \cup (\arcSetIn{\linePar} \cap \arcSetIn{u})\ \cup \notag \\
    & && \;  \cup (\arcSetChangeIn{\linePar} \cap \arcSetChangeIn{u}) \cup (\arcSetChangeOut{\linePar} \cap \arcSetChangeOut{u}) \notag \\
    & \flowVar{\arcPar}{\stationPar}{\periodPar{}} \leq \frequencyVar{\linePar}{\frequencyPar}{\periodPar{}} \cdot \frequencyPar \cdot \capacityPerLine{\linePar} \cdot \symmetryMultiplier  &&\forall \linePar \in \linepool,\ \stationPar \in \stationsOnLine{\linePar}\setminus \{\lastStation{\linePar}\},\ \arcPar \in (\arcSetIn{\linePar} \cap \arcSetIn{\stationPar}),\ \frequencyPar \in \frequencySet{\linePar},\ \periodPar{} \in \periodSet \label{MM-MP:useFrequencyInArcs}\\
    & \sum_{\stationPar \in \stationSet} \flowVar{\arcPar}{\stationPar}{\periodPar{}} \leq \frequencyVar{\linePar}{\frequencyPar}{\periodPar{}} \cdot \frequencyPar \cdot \capacityPerLine{\linePar} \cdot \symmetryMultiplier \quad &&\forall \linePar \in \linepool,\ \arcPar \in \arcSetChangeIn{\linePar},\ \frequencyPar \in \frequencySet{\linePar},\ \periodPar{} \in \periodSet \label{MM-MP:useFrequencyArcs}\\
    & \symmetryMultiplier \cdot \Bigl(\sum_{\linePar \in \linepool} \sum_{\frequencyPar \in \frequencySet{\linePar}} \costOfLine{\linePar} \cdot \frequencyPar \cdot \frequencyVar{\linePar}{\frequencyPar}{\periodPar{}} \Bigr) \leq \budget{\periodPar{}} \quad &&\forall \periodPar{} \in \periodSet \label{MM-MP:budgetConstraint} \\
    &\sum_{\linePar \in \linesLeavingStation{\stationPar}} \sum_{\periodPar{} \in \periodSet} \sum_{\frequencyPar \in \frequencySet{\linePar}} \periodLength{\periodPar{}} \cdot \frequencyPar \cdot \frequencyVar{\linePar}{\frequencyPar}{\periodPar{}} = \sum_{\linePar \in \linesEnteringStation{\stationPar}} \sum_{\periodPar{} \in \periodSet}  \!\!\!\!\!\!\! \! && \sum_{\frequencyPar \in \frequencySet{\linePar}} \periodLength{\periodPar{}} \cdot  \frequencyPar \cdot \frequencyVar{\linePar}{\frequencyPar}{\periodPar{}} \qquad \forall \stationPar \in \terminalStationSet \label{MM-MP:setEqualFrequencies} \\
    &\symmetryMultiplier \cdot \left(\stopVar{\linePar}{\stationPar}{\periodPar{j}} - \stopVar{\linePar}{\stationPar}{\periodPar{(j-1)}}\right) \leq \changeStopVar{\linePar}{\periodPar{j}}{\stationPar} \quad &&\forall \linePar \in \linepool,\ j \in \{2, \hdots,|\periodSet|\},\ \stationPar \in \stationsOnLine{\linePar} \setminus \{\firstStation{\linePar}, \lastStation{\linePar}\} \label{MM-MP:stopChange1}\\
    &\symmetryMultiplier \cdot \left(\stopVar{\linePar}{\stationPar}{\periodPar{(j-1)}} - \stopVar{\linePar}{\stationPar}{\periodPar{j}}\right) \leq \changeStopVar{\linePar}{\periodPar{j}}{\stationPar} \quad &&\forall \linePar \in \linepool,\ j \in \{2, \hdots,|\periodSet|\},\ \stationPar \in \stationsOnLine{\linePar} \setminus \{\firstStation{\linePar}, \lastStation{\linePar}\} \label{MM-MP:stopChange2}\\
    &\symmetryMultiplier \cdot \left(\frequencyVar{\linePar}{\frequencyPar}{\periodPar{j}} - \frequencyVar{\linePar}{\frequencyPar}{\periodPar{(j-1)}}\right) \leq \changeFreqVar{\linePar}{\periodPar{j}} \quad &&\forall \linePar \in \linepool,\ j \in \{2, \hdots,|\periodSet|\},\ \frequencyPar \in \frequencySet{\linePar} \label{MM-MP:freqChange1}\\
    &\symmetryMultiplier \cdot \left(\frequencyVar{\linePar}{\frequencyPar}{\periodPar{(j-1)}} - \frequencyVar{\linePar}{\frequencyPar}{\periodPar{j}}\right) \leq \changeFreqVar{\linePar}{\periodPar{j}} \quad &&\forall \linePar \in \linepool,\ j \in \{2, \hdots,|\periodSet|\},\ \frequencyPar \in \frequencySet{\linePar} \label{MM-MP:freqChange2}\\
    & \sum_{\linePar \in \linepool} \sum_{j=2}^{|\periodSet|} \Bigg( \sum_{\stationPar \in \stationsOnLine{\linePar} \setminus \{\firstStation{\linePar}, \lastStation{\linePar}\}} \changeStopVar{\linePar}{\periodPar{j}}{\stationPar} + \sum_{\frequencyPar \in \frequencySet{\linePar}} \!\!\!\!\!\!\!\!\! && \changeFreqVar{\linePar}{\periodPar{j}}\Bigg) + \slackVar = \varepsilon \label{MM-MP:maxChange} \\
    & \stopVar{\linePar}{\stationPar}{\periodPar{}} \in \{0,1\} \quad &&\forall \linePar \in \linepool,\ \stationPar \in \stationsOnLine{\linePar} \setminus \{\firstStation{\linePar}, \lastStation{\linePar}\},\ \frequencyPar \in \frequencySet{\linePar},\ \periodPar{} \in \periodSet \label{MM-MP:signStopVar} \\
    & \frequencyVar{\linePar}{\frequencyPar}{\periodPar{}} \in \{0,1\} \quad  &&\forall \linePar \in \linepool,\ \frequencyPar \in \frequencySet{\linePar},\ \periodPar{} \in \periodSet \label{MM-MP:signFrequencyVar} \\
    & \flowVar{\arcPar}{\stationPar}{\periodPar{}} \geq 0 \quad &&\forall \arcPar \in \arcSet{},\ \stationPar \in \stationSet,\ \periodPar{} \in \periodSet \label{MM-MP:signFlowVar} \\
    & \changeStopVar{\linePar}{\periodPar{j}}{\stationPar} \geq 0 \quad &&\forall \linePar \in \linepool,\ j \in \{2, \hdots,|\periodSet|\},\ \stationPar \in \stationsOnLine{\linePar} \setminus \{\firstStation{\linePar}, \lastStation{\linePar}\} \label{MM-MP:stopChangeSign} \\
    & \changeFreqVar{\linePar}{\periodPar{j}} \geq 0 \quad &&\forall \linePar \in \linepool,\ j \in \{2, \hdots,|\periodSet|\},\ \frequencyPar \in \frequencySet{\linePar} \label{MM-MP:freqChangeSign}
\end{align}

Constraint \eqref{MM-MP:flowConstraints} ensures that the passenger flow is conserved at all nodes. 
Next, \eqref{MM-MP:maxOneFrequency} arranges that a maximum of one frequency is chosen for each line.
Constraint \eqref{MM-MP:useDriveArc} ensures that the passenger flow on a Drive arc cannot exceed the capacity of the number of trains operated on this line. 
Therefore, if no frequency is selected, this line cannot be used to transport any passengers.
Constraint \eqref{MM-MP:useStopArc} ensures that when line \(\linePar\) will not stop at station \(\stationPar\), the Stop arc cannot be used, while \eqref{MM-MP:useSkipArc} ensures that when \(\linePar\) does stop at station \(\stationPar\), the Skip arc cannot be used.
Furthermore, constraint \eqref{MM-MP:useEnterOrLeaveLineArcs} assures that passengers can only enter or leave a train line when the train line stops at that station. 
To avoid creating quadratic constraints we use parameter \(\maxFreq{\linePar}\) in constraints \eqref{MM-MP:useStopArc}-\eqref{MM-MP:useEnterOrLeaveLineArcs} instead of \(\sum_{\frequencyPar \in \frequencySet{\linePar}} \frequencyPar \cdot \frequencyVar{\frequencyPar}{\linePar}{\periodPar{}}\) as was done in \eqref{MM-MP:useDriveArc}. 
Moreover, the capacities are multiplied by \(\symmetryMultiplier\), which means that the capacity is doubled when we consider symmetric lines.
This is needed as in that case these arcs are used for travelling in both directions.
The capacity restriction on Drive arcs by \eqref{MM-MP:useDriveArc} in combination with the flow conservation constraints \eqref{MM-MP:flowConstraints} ensure that the flow on the Skip arcs, Stop arcs, and arcs for entering or leaving a line will not exceed the capacity of the vehicles. 
Similarly, constraints \eqref{MM-MP:useFrequencyInArcs} and \eqref{MM-MP:useFrequencyArcs} make sure that an arc from a station's In or Change node to a train line can only be used if the corresponding frequency is selected for the train line.
Constraint \eqref{MM-MP:budgetConstraint} ensures that the cost of the selected line plan does not exceed the given budget. 
Here we define the cost of a period's line plan as the total number of kilometres driven by the trains in one hour. 
These train kilometres are calculated by multiplying the length of each train line in one direction \((\costOfLine{\linePar})\) by the chosen frequency for that line, and then summing over all lines.
If we consider symmetric lines, each line is driven in both direction.
Therefore, in that case we multiply the one-way cost of each line by two via the symmetry parameter \(\symmetryMultiplier\).
Constraint \eqref{MM-MP:setEqualFrequencies} enforces that for each terminal station \(\stationPar \in \terminalStationSet\), the total frequency of lines leaving this station should be equal to the total frequency of the lines entering this station.
By balancing the number of trains leaving and entering the terminal station, it is more likely that in a later planning stage a rolling stock plan can be found.
Note that in the case of symmetric lines, this constraint can be removed as it will always be satisfied. 
The reason for this is that symmetric lines are operated in both directions, and hence the same number of trains will depart and arrive at each terminal station on the line.

Constraints \eqref{MM-MP:stopChange1}-\eqref{MM-MP:maxChange} calculate and limit the number of adjustments made to the line plan throughout the day.
We limit the number of adjustments, as making changes to the line plan during the day has a cost for both the passengers and the railway undertaking.
Passengers value regularity in the timetable, because it makes it easy to use and remember the schedule.
\cite{wardman2004ConsumerBenefitsDemand} and \cite{johnson2006ForecastingAppraisingImpact} show that a regular timetable increases the passenger demand. 
Furthermore, every change to the line plan affects the timetable and possibly the rolling stock schedule. 
Since there is still a lot of manual labour involved in the public transport planning process, having many adjustments would mean a significant increase in cost for the railway undertaking.

In this work, we measure the dissimilarity between two successive line plans by counting the number of changes made to the stopping patterns and frequencies of the lines. 
Previous work of \cite{duran-micco2022DesigningBusLine} and \cite{schiewe2023LinePlanningDifferent} on line planning for bus and metro also provide metrics for measuring the (dis)similarity between different line plans.
However, these metrics are mainly based on the common edges in both line plans, which makes them less suitable for appropriately evaluating the changes to the stopping pattern for railways.
When a train would skip a stop, its route through the network would likely not change.
Hence, we count the number of changes to the stopping pattern and frequencies instead.

Constraints \eqref{MM-MP:stopChange1} and \eqref{MM-MP:stopChange2} increase the lower bound of the change stop variable \(\changeStopVar{\linePar}{\periodPar{j}}{\stationPar}\) to \(\symmetryMultiplier\) if the value of the corresponding stop variables is not the same in two consecutive periods.
Similarly, constraints \eqref{MM-MP:freqChange1} and \eqref{MM-MP:freqChange2} set the lower bound of the change frequency variable \(\changeFreqVar{\linePar}{\periodPar{j}}\) to \(\symmetryMultiplier\) if a different frequency is chosen in two consecutive periods.
Note that each adjustment of the stopping pattern or frequency is counted once when asymmetric lines are considered.
For symmetric lines each adjustment counts twice, as in that case the changes are made in both directions.
Constraint \eqref{MM-MP:maxChange} sums all change variables and together with a slack variable \(\slackVar\) this sum should be equal to a maximum number of changes \(\varepsilon\).
Lastly, constraints \eqref{MM-MP:signStopVar} and \eqref{MM-MP:signFrequencyVar} denote that the stop variables and frequency variables are binary, and constraints \eqref{MM-MP:signFlowVar}-\eqref{MM-MP:freqChangeSign} indicate that the flow, change stop, and change frequency variables should be non-negative.

\subsubsection{Objective}
The model's aim is to create a line plan that minimises the passengers' GJT, which can be calculated using \eqref{MM-MP:obj_GJT}.
Furthermore, as we want to find solutions that do not have more adjustments to the line plan than necessary, we add a cost term to maximise the value of the slack variable multiplied by a small coefficient \(\delta\). 
This gives us the following optimisation model:
\begin{equation} \label{model-multiperiod}
    \begin{aligned}
        \text{Minimise} \quad &\sum_{\periodPar{} \in \periodSet} \sum_{\arcPar \in \arcSet{}} \sum_{\stationPar \in \stationSet} \periodLength{\periodPar{}} \arcCost{\arcPar} \flowVar{\arcPar}{\stationPar}{\periodPar{}} - \delta \cdot \slackVar\\
        \text{s.t.} \quad &\text{Constraints } \eqref{MM-MP:flowConstraints}-\eqref{MM-MP:freqChangeSign}
    \end{aligned} \tag{MP-LPP}
\end{equation}
By solving this model for varying bounds on the number of allowed changes \((\varepsilon)\), we can generate a set of Pareto optimal line plans.
This provides a decision maker at the railway undertaking the opportunity to investigate multiple plans and select the one that is most to their liking. 
This so-called \eps-constraint method is a popular method for solving multi-objective integer linear programming models \citep{mavrotas2013ImprovedVersionAugmented} and has been applied several times in the railway planning domain (see e.g., \cite{binder2017MultiobjectiveRailwayTimetable}, \cite{yan2019MultiobjectivePeriodicRailway}, \cite{stoilova2020IntegratedMultiCriteriaMultiObjective}, \cite{ning2022BiobjectiveOptimizationModel},\cite{yang2024IntegratingStopPlanning}).
In this work, We use the AUGMECON2 algorithm by \cite{mavrotas2013ImprovedVersionAugmented} to approximate the Pareto optimal solutions.

\section{Case study description}\label{sec:case_study}
\begin{figure}
    \centering
    \includegraphics[width=\textwidth]{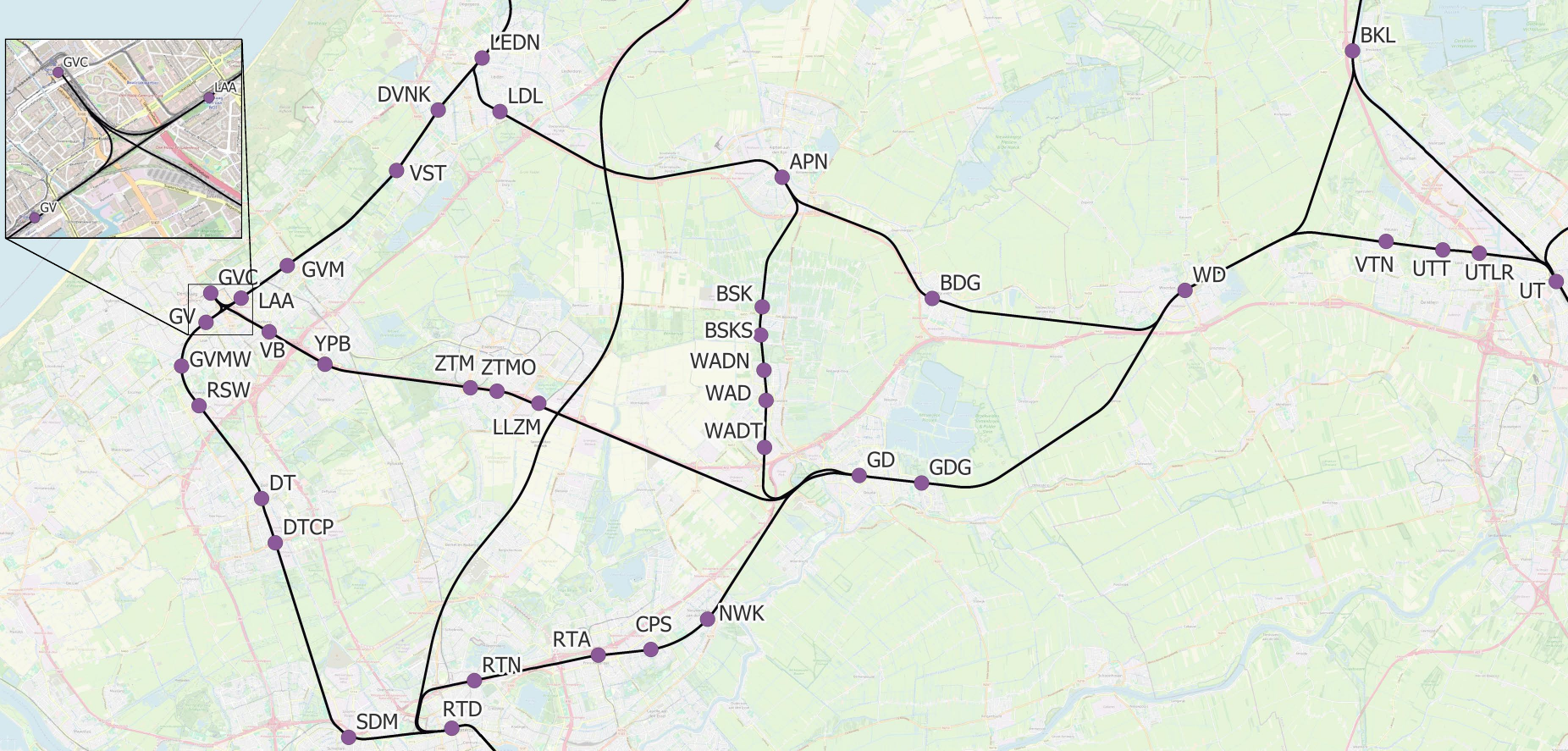}
    \caption{Part of the Dutch railway network between the cities Leiden (Ledn), The Hague (Gvc), Rotterdam (Rtd), and Utrecht (Ut).}
    \label{fig:network_Rtd-Gvc-Ledn-Ut}
\end{figure}

We test the model described in Section \ref{sec:mathematicalModel} on a case study based on part of the railway network in the Netherlands.
We consider the network between the cities Leiden, The Hague, Rotterdam, and Utrecht.
Rotterdam, The Hague, and Utrecht are the second, third, and fourth-largest municipalities of the Netherlands with populations between 673,000 and 377,000 \citep{cbs2025InwonersGemeente}. 
Leiden is a smaller city with around 131,000 inhabitants, which makes it the third-largest city in the province South-Holland, after Rotterdam and The Hague.
This network was selected as it has a lot of (directional) peak demand and different types of stations, including main stations of large cities and stations that are mostly interesting for commuters.
Therefore, it is likely that different line plans are optimal in different periods, and hence that adjustments will be made during the day.
A graphical representation of this network is depicted in Figure \ref{fig:network_Rtd-Gvc-Ledn-Ut}.

In Table \ref{tab:stationAbrreviations} in Appendix \ref{ap:caseStudy} an overview is given of the abbreviations used in Figure \ref{fig:network_Rtd-Gvc-Ledn-Ut} and the corresponding station names. 
This table also denotes if a station is a terminal station (E) and if transferring between trains is possible at the station (T).
We have tested different sets of transfer stations and these tests showed that when transfers are possible at all stations, the same line plan was found as when transfers are only allowed at the larger stations.
However, the solver needed either more time to solve the instance or had a larger optimality gap when transfers were possible at all stations. 
Therefore, we choose to only allow for transfers at the large stations in the network. 

\begin{table}[t!]
    \centering
    \small
    \caption{Set of candidate lines for the case study with symmetric lines.}
    \label{tab:linepool_UtRtdGvcLedn}
    \begin{tabular}{llccc}
    \hline \Tstrut
    Name \((\linePar)\) & Stations on route (\(\stationsOnLine{\linePar}\)) & \(\frequencySet{\linePar}\) & \(\costOfLine{\linePar}\) & \(\capacityPerLine{\linePar}\)  \\
    &&  [trains/hr/dir] & [km] & [seats/train]\Bstrut \\ \hline \Tstrut
    Gd-Apn    & Gd,Wadt,Wad,Wadn,Bsks,Bsk,       & \{1,2,4\}     & 17.5  & 314 \\
    & Apn && \\
    Gvc-Gdgo  & Gvc,Vb,Gvy,Ztm,Ztmo,Llzm,        & \{1,2,4,6\}   & 30.5  & 816 \\
    & Gd,Gdgo && \\
    Gvc-Ledn  & Gvc,Laa,Gvm,Vst,Dvnk,Ledn        & \{1,2,4,6\}   & 15.6  & 816 \\
    Gvc-Ut    & Gvc,Vb,Gvy,Ztm,Ztmo,Llzm,        & \{1,2,4,6\}   & 60.1  & 816 \\
    & Gd,Gdgo,Wd,Vtn,Utt,Utlr,Ut && \\
    Ledn-Ut   & Ledn,Ldl,Apn,Bdg,Wd,Vtn,         & \{1,2,4\}     & 49.1  & 816 \\
    & Utt,Utlr,Ut && \\
    Rtd-Bkl   & Rtd,Rtn,Rta,Cps,Nwk,Gd,Gdgo,     & \{1,2,4,6\}   & 52.8  & 628 \\
    & Wd,Bkl && \\
    Rtd-Gdgo  & Rtd,Rtn,Rta,Cps,Nwk,Gd,Gdgo      & \{1,2,4,6\}   & 26.0  & 628 \\
    Rtd-Gvc   & Rtd,Sdm,Dtcp,Dt,Rsw,Gvmw,        & \{1,2,4,6\}   & 23.2  & 816 \\
    & Gv,Gvc & & & \\
    Rtd-Ledn  & Rtd,Sdm,Dtcp,Dt,Rsw,Gvmw,        & \{1,2,4,6\}   & 37.6  & 816 \\
    & Gv,Laa,Gvm,Vst,Dvnk,Ledn & & & \\
    Rtd-Ut    & Rtd,Rtn,Rta,Cps,Nwk,Gd,Gdgo,     & \{1,2,4,6\}   & 55.6  & 628 \\
    & Wd,Vtn,Utt,Utlr,Ut && \Bstrut \\  \hline
    \end{tabular}
\end{table}
The set of candidate lines used in the case study is given in Table \ref{tab:linepool_UtRtdGvcLedn}. 
The first column displays the names of the candidate lines, which consist of the first station and last station.
A subroute with an earlier terminal is considered a different route, since terminals cannot be skipped in the stop pattern.
As we only consider one candidate line per route in this case study, the names are unique.
Tests with multiple candidate lines per route indicated that the model prefers one line per route with a high frequency over two lines with different stopping patterns and lower frequencies.
Since we have a smaller solution space with one candidate line per route, and therefore can find better solutions in the same running time, we have chosen to include only one candidate line per route.
If more lines per route are considered, an index can be added to create unique line names.
The second column of Table \ref{tab:linepool_UtRtdGvcLedn} denotes the sequence of stations that a train passes in the network. 
This sequence determines the routing through the public transport network. 
Next, the third column displays the set of frequencies that a line is allowed to take. 
The routes of lines between Gouda (Gd) and Alphen a/d Rijn (Apn) and between Leiden (Ledn) and Utrecht (Ut) partially go over single-track sections in the network.
Therefore, the maximum frequency of these lines is set to 4 per hour instead of 6.
The second to last column displays the (one-way) length of the line in kilometres, which is taken as proxy for the cost of the line and the last column shows the seat capacity of trains on this line.
This capacity is based on the largest train that can stop at all stations on a line's route. 
Note that Table \ref{tab:linepool_UtRtdGvcLedn} only shows the symmetric lines.
For the case with asymmetric lines the reverse of these lines are also added to the set of candidate lines with the same frequency set, cost and capacity.

In this case study, 3 periods are considered: the morning peak, the midday off-peak, and the afternoon peak.
By considering these periods, we have a large variation in demand.
Passenger volumes are much larger in the peak periods than in the off-peak period and by looking at both the morning and afternoon peak, we also see a change in the direction of the passenger flows.
For each period, a 1-hour OD matrix that represents the period's demand is needed as input for the model.
In previous work of the authors \citep{vanderknaap2024ClusteringRailwayPassenger}, a method was developed to determine periods in which the railway demand is homogeneous. 
To apply the results of \cite{vanderknaap2024ClusteringRailwayPassenger}, we use the same demand data set based on the realised demand of NS on Tuesdays in 2019.
NS obtains this demand data from a smart card system with check-in and check-out gates.
More details about the data set can be found in \cite{vanderknaap2024ClusteringRailwayPassenger} and details about this smart card system are provided in \cite{vanoort2015ShorttermPredictionRidership}.

According to the study, the morning and afternoon peaks both have three different types of demand: start of the peak, hyper peak, and end of the peak.
We choose to use the demand during the hyper peaks (in morning from 7:30-8:30 and in the afternoon from 16:30-18:00) to represent the peak demand.
By using the demand from the busiest time during the peak, it is more likely that we find a line plan that can serve all the demand. 
Furthermore, \cite{vanderknaap2024ClusteringRailwayPassenger} found the demand characterising the midday off-peak to be the demand between 9:30 and 15:00.
As the selected demand periods are not connected, we have to decide what schedule should be used in the periods 8:30-9:30 and 15:00-16:30.
According to \cite{vanderknaap2024ClusteringRailwayPassenger}, the demand during 8:30-9:30 is closer to the morning hyper peak than to the midday off-peak demand. 
Therefore, we propose to use the schedule found for the morning peak until the off-peak starts at 9:30.
Similarly, the demand between 15:00 and 16:30 is closer to the demand in the afternoon hyper peak than the demand of the off-peak, so the schedule for the afternoon peak can be used from 15:00 onwards. 
To create the 1-hour OD matrix for each period, we sum the half hour OD matrices in the data set of which the period consists and then divide this by the number of hours in the period. 
In the remainder of this paper we will refer to the morning hyper peak period with \mop, the midday off-peak period with \mo, and the afternoon hyper peak period with \ap.
The period lengths are \(\periodLength{\text{\mop}}= 1.5\), \(\periodLength{\text{\mo}}=5.5\), and \(\periodLength{\text{\ap}}=1.5\). 
Note that we increased the length of \mop\ by 0.5, such that \mop\ and \ap\ have the same length. 
This allows the model to create line plans that have imbalanced frequencies within periods but balanced frequencies over the whole day.
As the passenger volumes in the peak hours are much higher than the volumes during the off-peak, this is also reflected in the budgets for the line plans.
We set \(\budget{\text{\mop}}=\budget{\text{\ap}}=3108\), based on the cost of the line plan operated on this part of the network in May 2023, and we use half of this budget for the off-peak \((\budget{\text{\mo}}=1554)\) due to the significantly lower demand.

Lastly, we need to determine the coefficients used in the objective of \eqref{model-multiperiod}.
To calculate the GJT, we need to give a cost for each arc in the C\&G network.
For Drive arcs, the arc cost is the driving time between two stations, which is calculated using the length of each line section and the average speed on that section.
For all other arcs, the costs are given in Table \ref{tab:arcCosts_caseStudy}.
\begin{table}[t]
    \centering
    \caption{Used costs of arcs in the C\&G-network. Costs are given in generalised journey time minutes.}
    \label{tab:arcCosts_caseStudy}
    \begin{tabular}{llll}
        \hline \Tstrut
        Arc type   & Cost of arc [min] & Arc type     & Cost of arc [min] \Bstrut \\ \hline \Tstrut
        In-F1      & 55.85             & In-change-F1 & 50.25             \\
        In-F2      & 31.85             & In-change-F2 & 28.85             \\
        In-F4      & 17.35             & In-change-F4 & 19.95             \\
        In-F6      & 10.85             & In-change-F6 & 15.55             \\
        Out        & 0.7               & Skip         & 0                 \\
        Out-change & 0.7               & Stop         & 3.55             \Bstrut \\ \hline
    \end{tabular}
\end{table}
\begin{figure}[t]
    \centering
    \includegraphics[width=0.9\textwidth]{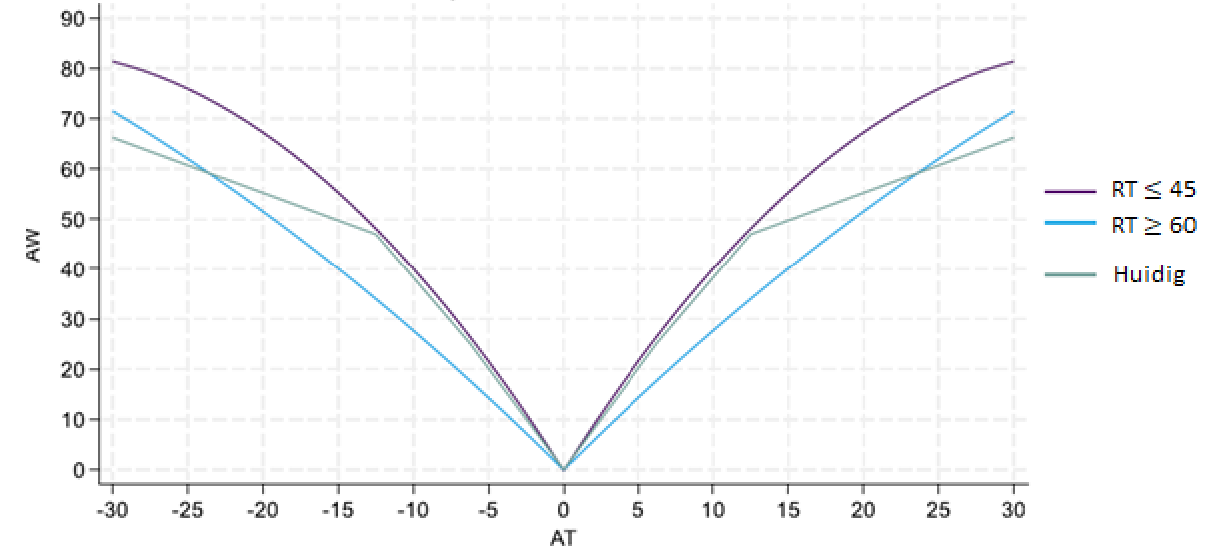}
    \caption{Entry resistance curve from \cite{guis2023ReizigerKiestKorter}. Vertical axis denotes the entry resistance (in minutes), while the horizontal axis denotes the passenger's required adaptation from the desired departure time (in minutes). Values of the purple curve for travel times (denoted by RT in the legend) up to 45 minutes is used.}
    \label{fig:entry_resistance}
\end{figure}
The cost of the In and In-change arcs are dependent on the frequency of the line.
For example, in Table \ref{tab:arcCosts_caseStudy} In-F2 denotes the cost of the in arc if the frequency of the line is two. 
The values in Table \ref{tab:arcCosts_caseStudy} of the In and In-change arcs include an acceleration penalty of 0.85 minutes to account for the increased travel time due to acceleration after the stop.
Similarly, the costs of the Out and Out-change arcs include a deceleration time penalty of 0.7 minutes, to account for the extra travel time needed due to braking.  
The value of the Stop arc includes a dwell time of 2 minutes and both the acceleration and deceleration penalties. 
The cost of the In and In-change arcs are based on the work of \cite{guis2023ReizigerKiestKorter} and \cite{debruyn2023ReizigerKiestKorter}, respectively.
Both papers used a stated preference experiment to respectively estimate the entry and transfer resistance in terms of GJT for train passengers in the Netherlands. 
\cite{guis2023ReizigerKiestKorter} provide an entry resistance curve, which is shown in Figure \ref{fig:entry_resistance}.
For each adaptation that a passenger has to make from the preferred departure time, this entry resistance curve provides the generalised journey time that the passenger experiences because of this change.
Since the timetable is not yet known in this stage, and hence neither the adaptation time for the passengers, we assume that the passengers want to depart according to a uniform distribution.
Furthermore, we assume that the trains on each line are equally spaced over the hour.
Given these assumptions, if a train has a frequency of two times per hour, the maximal adaptation time would be 15 minutes (earlier or later than the preferred departure time).
If the passengers' preferred departure times are uniformly distributed, then the average adaptation time would be 7.5 minutes for a frequency of two. 
Using the entry resistance curve provided by \cite{guis2023ReizigerKiestKorter}, we find that the cost of the corresponding in arc is equal to 31 GJT minutes, to which the previously discussed acceleration penalty of 0.85 is added. 
Similarly, the values of the other In arcs in Table \ref{tab:arcCosts_caseStudy} are derived from this curve. 
To calculate the GJT of a transfer, \cite{debruyn2023ReizigerKiestKorter} derive a function that takes into account the type of transfer (e.g., cross-platform), the waiting time for the connecting train, and the time you have to wait for the next train if you miss the transfer.
We use this function to calculate the values of the In-change arcs displayed in Table \ref{tab:arcCosts_caseStudy}, under the assumption that the transfer will be cross-platform, the waiting time is half the headway, and the extra waiting time in case of a missed transfer is 60 divided by the frequency (assuming regular services).
Lastly, according to \cite{mavrotas2013ImprovedVersionAugmented} the \(\delta\) in the objective of \ref{model-multiperiod} should take a very small value (usually \(10^{-3}\) to \(10^{-6}\)).
\(\delta\) should be small such that it does not influence the main objective (in our case minimising the GJT).
Moreover, the value should be large enough such that we find the solution with the least amount of changes to the line plan, when multiple solutions with the same GJT exist.
For the case study we set \(\delta = 10^{-3}\).
The results of the case study are provided in the next section.

\section{Results}\label{sec:results}
In this section we demonstrate the model described in Section \ref{sec:mathematicalModel} by applying it to a case study using real data of part of the Dutch railway network as described in Section \ref{sec:case_study}. 
The model was coded in Python 3.10.11 and solved with Gurobi Optimizer 10.0.1. 
All experiments were carried out on a laptop with Intel\textsuperscript{\textregistered} Core\texttrademark \ i7-1185G7 @ 3.00GHz and 16 GB RAM.

\subsection{General insights}\label{subsec:gen_results}
Using the \eps-constraint method, 10 different multi-period line plans with symmetric lines and 13 plans with asymmetric lines are generated by solving \eqref{model-multiperiod} with different bounds \eps\ on the number of line plan adjustments in constraint \eqref{MM-MP:maxChange}.
Table \ref{tab:results_Ut_case} shows for each plan the type of lines used (symmetric or asymmetric), the total number of adjustments made to the line plan and how many of these are adjustments of the frequencies and adjustments of the stopping pattern.
We count an adjustment to the line plan if the stopping pattern at a station or the frequency of a line is different relative to the previous period. 
Hence, in this case the total number of adjustments to the line plan is the number of adjustments in period 2 relative to period 1 plus the adjustments in period 3 relative to period 2.
The last three columns of Table \ref{tab:results_Ut_case} show the total GJT, the \%-change in total GJT compared the GJT of our reference scenario (denoted in bold), and the optimality gap of the solution.
For calculating the optimality gap, we have used two options. 
The first option is the optimality gap reported by Gurobi after reaching the time limit. 
These gaps are denoted with a * in Table \ref{tab:results_Ut_case}.
As the gaps reported by Gurobi were quite large, we have also calculated a different lower bound by optimising the line plan for each period separately.
Since the constraints concerning multiple periods (i.e., \eqref{MM-MP:setEqualFrequencies}-\eqref{MM-MP:maxChange}) are not taken into account during these optimisation problems, they provide a lower bound for each period.

As reference scenario, we have selected the scenario with symmetric lines and 20 adjustments.
In 2023, only symmetric lines were operated on the network, which is why we take the symmetric lines as reference.
Furthermore, we need to allow for some changes in the reference scenario as in the midday off-peak (\mo) only half the amount of train kilometres are allowed compared to the morning peak (\mop) and afternoon peak (\ap).
Therefore, several adjustments of the frequency are needed to have line plans in \mop\ and \ap\ that (approximately) use the allowed budget.
The scenario with 20 adjustments is the first scenario that can use at least 98.5\% of the budget  during each period, which is why we use this as the reference scenario.
We can also see this need for changing the frequency in the third and fourth column of Table \ref{tab:results_Ut_case}.
The first six scenarios with symmetric lines only change the frequency and the first five scenarios with asymmetric lines also mainly focus on changing the frequency.
After that, the focus shifts to making changes to the stopping pattern.

\begin{table}[t!]
    \centering
    \caption{Results of different line plans created by the \eps-constraint method with a time limit of 2 hours (7200 seconds). The base scenario is the scenario with symmetric lines and 20 adjustments total (denoted in bold). Optimality gaps with a * were reported by Gurobi.}
    \label{tab:results_Ut_case}
    \begin{tabular}{lcccccc}
        \hline
        Type of     & \multicolumn{3}{c}{Number of adjustments}& Total         & \% difference & Optimality      \\ \cline{2-4}
        lines       & Total      & Freq.      & Stopping  & GJT                & in GJT w.r.t. & gap       \\
                    &            &            & pattern   & [min]              & base scenario &           \\ \hline
        Symmetric         & 0          & 0          & 0         & 6,744,095          & 19.74\%   & 4.8\%*        \\
        Symmetric         & 4          & 4          & 0         & 6,150,780          & 9.20\%    & 16.5\%        \\
        Symmetric         & 10         & 10         & 0         & 5,857,551          & 4.00\%    & 10.9\%        \\
        Symmetric         & 14         & 14         & 0         & 5,730,032          & 1.73\%    & 8.5\%         \\
        \textbf{Symmetric}& \textbf{20}& \textbf{20}& \textbf{0}& \textbf{5,632,353} & \textbf{-}& \textbf{6.7\%}\\
        Symmetric         & 24         & 24         & 0         & 5,576,424          & -0.99\%   & 5.6\%         \\
        Symmetric         & 30         & 24         & 6         & 5,529,375          & -1.83\%   & 4.7\%         \\
        Symmetric         & 32         & 24         & 8         & 5,527,530          & -1.86\%   & 4.7\%         \\
        Symmetric         & 36         & 24         & 12        & 5,527,296          & -1.87\%   & 4.7\%         \\
        Symmetric         & 40         & 24         & 16        & 5,523,326          & -1.94\%   & 4.6\%         \\ \hline
        Asymmetric        & 0          & 0          & 0         & 6,574,592          & 16.73\%   & 26.1\%*       \\
        Asymmetric        & 5          & 4          & 1         & 6,119,410          & 8.65\%    & 26.5\%*       \\
        Asymmetric        & 10         & 10         & 0         & 5,804,499          & 3.06\%    & 22.6\%        \\
        Asymmetric        & 15         & 14         & 1         & 5,602,865          & -0.52\%   & 18.0\%        \\
        Asymmetric        & 20         & 20         & 0         & 5,532,751          & -1.77\%   & 16.6\%        \\
        Asymmetric        & 30         & 24         & 6         & 5,519,459          & -2.00\%   & 16.3\%        \\
        Asymmetric        & 35         & 24         & 11        & 5,501,986          & -2.31\%   & 15.9\%        \\
        Asymmetric        & 40         & 24         & 16        & 5,440,503          & -3.41\%   & 16.3\%        \\
        Asymmetric        & 44         & 22         & 22        & 5,440,403          & -3.41\%   & 16.1\%        \\
        Asymmetric        & 55         & 22         & 33        & 5,419,928          & -3.77\%   & 15.7\%        \\
        Asymmetric        & 60         & 22         & 38        & 5,419,109          & -3.79\%   & 15.5\%        \\
        Asymmetric        & 65         & 22         & 43        & 5,403,432          & -4.06\%   & 15.3\%        \\
        Asymmetric        & 70         & 26         & 44        & 5,392,625          & -4.26\%   & 15.2\%        \\ \hline
    \end{tabular}
\end{table}

In Table \ref{tab:results_Ut_case}, the reported line plans with the highest number of adjustments are the line plans that were found when no restrictions were imposed on the number of adjustments. 
For the case with symmetric lines, it is optimal to make 40 adjustments to the line plan, while for the case with asymmetric lines 70 adjustments are optimal.
When we compare the GJT of these line plans with the GJT of the reference line plan, we see a maximal reduction of 1.94\% in the GJT when considering symmetric lines and a maximal reduction of 4.26\% if we look at the asymmetric lines.  
Moreover, if we compare the total GJT of the scenarios with asymmetric lines for 70 and 20 adjustments, we see an improvement of 2.53\%.
Therefore, we can indeed say that the model can create a line plan with adjusted stops that better matches the demand (in terms of providing the lowest GJT possible), than the line plan that only makes minimal changes to the frequencies to cut the service in the off-peak.
Note that a 4.26\% reduction in GJT can have a significant impact on a railway undertaking's revenue.
According to a meta-analysis of time elasticities of travel demand in the United Kingdom, the GJT elasticity of trains is around -0.81 \citep{wardman2012ReviewMetaanalysisTime}.
This means that a 4.26\% reduction in GJT would increase the revenue by 3.45\%.
As the NS' revenue in 2023 in the Netherlands was \texteuro 2877 million \citep{NS2023report}, a 3.45\% increase would be an increase of \texteuro 99.3 million.
In the following two subsections we look more closely into the effect of several design choices.

\subsection{Impact of asymmetric lines}
The first design choice that we want to address is the choice for either asymmetric or symmetric lines. 
\begin{table}[b]
\centering
\caption{Total GJT for line plans with symmetric and asymmetric lines with the same number of adjustments.}
\label{tab:sym_VS_asym_results}
\begin{tabular}{cccc}
\hline
    Number of   & \multicolumn{2}{l}{Total GJT [min]}      &               \\ \cline{2-3}
adjustments & Symmetric lines & Asymmetric lines & \%-difference \\ \hline
0           & 6,744,095         & 6,574,592          & -2.5\%        \\
10          & 5,857,551         & 5,804,499          & -0.9\%        \\
20          & 5,632,353         & 5,532,751          & -1.8\%        \\
30          & 5,529,375         & 5,519,459          & -0.2\%        \\
40          & 5,523,326         & 5,440,503          & -1.5\%        \\ \hline
\end{tabular}
\end{table}
In Table \ref{tab:sym_VS_asym_results}, we compare the total GJT of the line plans with symmetric lines to the GJT of the line plans with asymmetric lines with the same number of adjustments.
In this table, the first column denotes the number of adjustments that the line plans have, the second (third) column denotes the total GJT of the line plans with the symmetric (asymmetric) lines, and the last column denotes the percentage difference between the symmetric and asymmetric results.
This last column shows that the line plans with asymmetric lines have a GJT that is between 0.2\% and 2.5\% lower than the GJT of line plans with symmetric lines and the same amount of adjustments.
So if we allow for asymmetric lines, we can create line plans with lower total GJT given the same budget and allowed number of changes.

In addition to achieving a lower GJT with the same budget and number of changes, considering asymmetric lines also allows for more opportunities to lower the GJT through changing the line plan.
For the scenarios with asymmetric lines 70 adjustments were found, while for the symmetric lines a line plan with only 40 changes was created when no restriction was given for the number of changes.
One reason for this is that in the asymmetric case, adjustments can also be made in only one direction without the need to mirror these changes in the opposite direction.
For example, if a station is skipped in the peak direction, then this station can still be served in the opposite direction to provide some service to these passengers.
This approach saves GJT for passengers moving in the peak direction while still delivering service to the skipped station's passengers.

It is not surprising that the model can find a line plan with a lower GJT when asymmetric lines are allowed.
After all, all the solutions with symmetric lines are also allowed in the asymmetric case, so the GJT should be at least as good as in the symmetric case.
However, due to the larger solution space, the performance in terms of optimality gap are also worse for the cases with asymmetric lines than for the cases with symmetric lines. 
As can be seen in Table \ref{tab:results_Ut_case}, the optimality gaps for the line plans with symmetric lines fall between 4.6 and 16.5\%, while those for the line plans with asymmetric lines fall between 15.2 and 26.5\%.
So although we can find line plans that have lower GJT and are better adjustable to peak demand when we allow for asymmetric lines, this also comes with the cost of larger optimality gaps and potentially higher required running times. 
This could be problematic when good solutions need to be found in a relatively short time window.

\subsection{Impact of frequency and stopping pattern adjustments}
\begin{figure}[b!]
    \centering
    \includegraphics[width=0.8\textwidth]{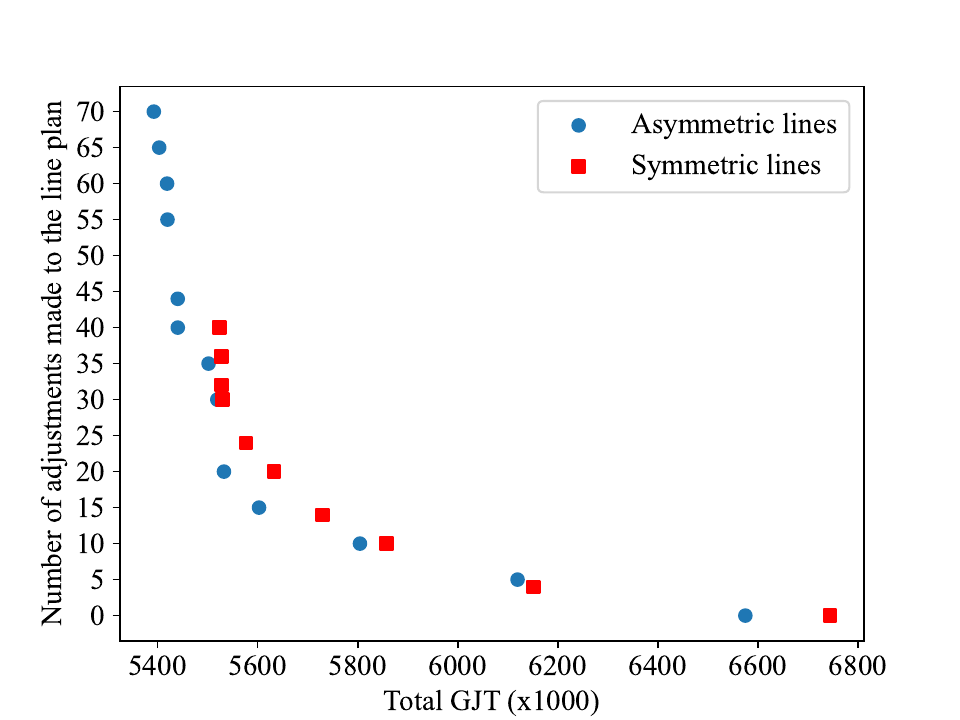}
    \caption{The approximated Pareto frontiers of the scenarios with symmetric lines and those with asymmetric lines}
    \label{fig:pareto_frontier}
\end{figure}

When looking at how much GJT an adjustment can save, we see that earlier adjustments yield larger GJT savings compared to later ones. 
This phenomenon is illustrated in Figure \ref{fig:pareto_frontier}, where the two approximated Pareto frontiers of the asymmetric and symmetric results are plotted. 
In this figure, the horizontal axis represents the total GJT and the vertical axis shows the number of adjustments made to the line plan.
Furthermore, the results for symmetric lines are marked by red squares, while asymmetric results are indicated by blue circles.
Figure \ref{fig:pareto_frontier} clearly shows small changes in the GJT when many adjustments are used and large changes in GJT when adjustments are few. 
For example, the case with asymmetric lines going from 0 to 5 changes, reduces the GJT by 6.9\%, while going from 40 to 44 changes only reduces the GJT by 0.002\%.
Given these diminishing GJT savings and the previously mentioned costs of adjustments for passengers and the railway undertaking, practitioners should carefully consider what the optimal number of changes is for their line plan.

Looking at the third and fourth column in Table \ref{tab:results_Ut_case}, we see that the first 24 changes made to the line plan are mainly focused on the frequencies, and the additional changes after that are mainly focussed on the stopping pattern. 
Hence, the frequency adjustments apparently have a larger impact on the GJT than stopping pattern adjustments.
This large impact can be explained, as the waiting time for boarding the first train and the transfer time are both frequency based. 
Since all passengers incur the cost of waiting to board the first train, all passengers that want to use a certain line will benefit if the frequency is increased. 
On the other hand, if a change is made to the stopping pattern, then some passengers benefit (e.g., because they have a lower in-vehicle time due to a skipped stop) while others incur higher costs (e.g., because they have to take a detour as their station is no longer served in their desired direction). 

\subsubsection{Line plan adjustments with symmetric lines}
\begin{figure}
    \centering
    \includegraphics[width=\textwidth]{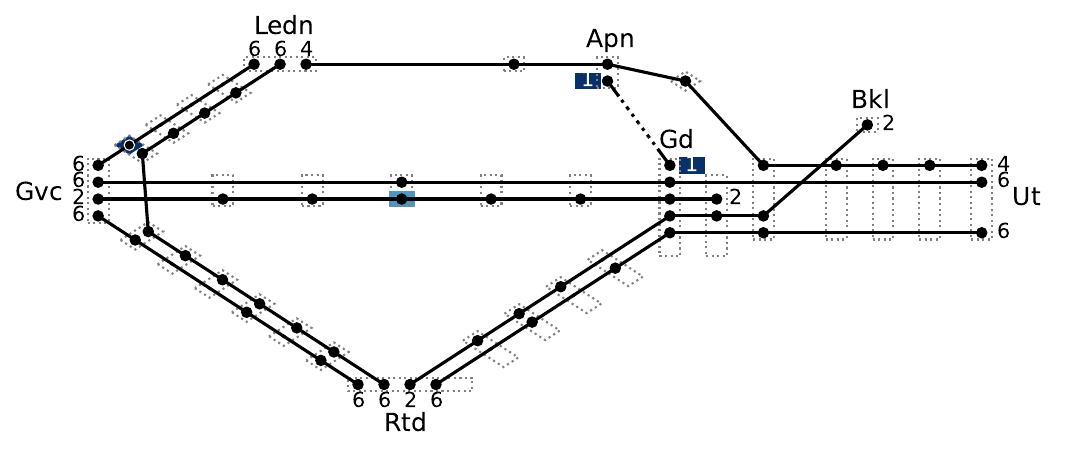}
    \includegraphics[width=\textwidth]{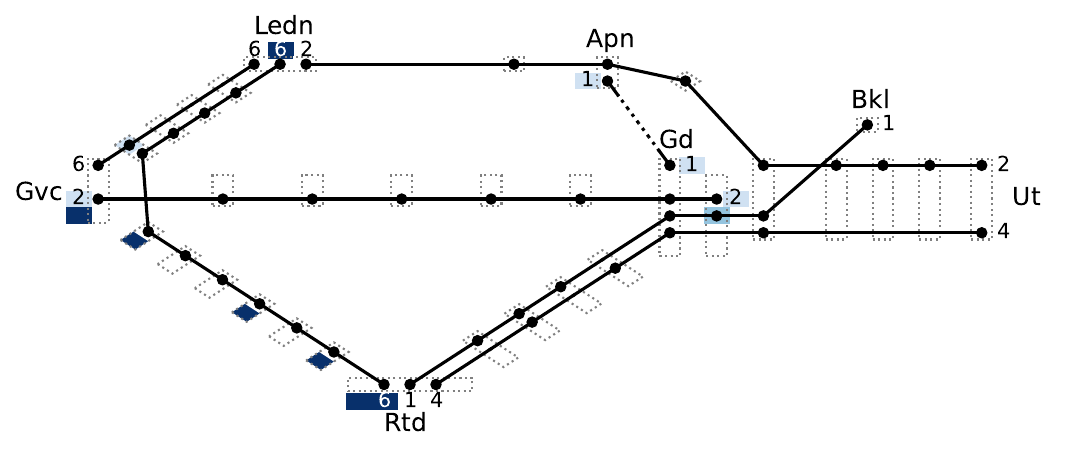}
    \includegraphics[width=\textwidth]{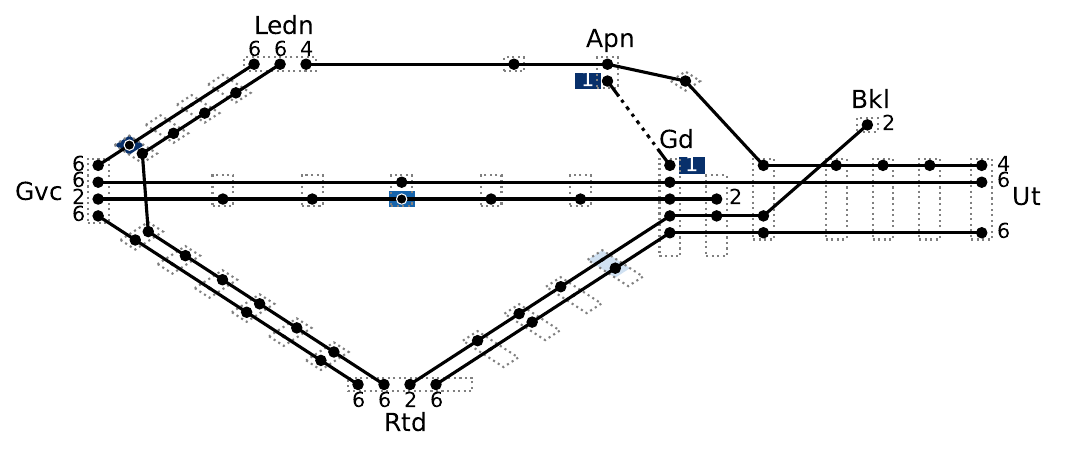}
    \caption{Line plan in \mop, \mo, and \ap\ (from top to bottom) with symmetric lines and 20 changes. Darker colours denote that the stopping pattern or frequencies were changed more often in scenarios with more changes. The line between Apn and Gd stops at all intermediate stations in all scenarios.}
    \label{fig:heatmap_sym}
\end{figure}

In Figure \ref{fig:heatmap_sym} we display the reference line plan with symmetric lines and 20 changes. 
In this figure, a circle denotes that a line stops at a station.
The frequencies of the lines are denoted at their terminal stations. 
Furthermore, we have compared this reference line plan to all line plans with symmetric lines and more than 20 adjustments. 
In Figure \ref{fig:heatmap_sym} we have also marked where the line plans are changed compared to the reference line plan.
The colour denotes how often the line plan is changed at that location, where white denotes that no changes were made, and dark blue indicates that the change was made in many scenarios.
The colour behind the frequency indicates how often this was changed in other scenarios.
Similarly, for the stopping pattern: if the background of a dot is coloured, then this stop is changed to a skip, while a coloured background of a skip means that the line stops at this station in other scenarios. 

When we compare the frequencies over the periods of the reference line plan, we see that the frequencies of most lines are reduced by one step (e.g., from 6 to 4 or from 4 to 2) in the off-peak, and increased again for the afternoon peak.
However, there are also two lines (between Rotterdam Centraal (Rtd) and Den Haag Centraal (Gvc) and between Gvc and Ut) that completely disappear in the off-peak.
Note that the services that disappear in the off-peak are the fast `Intercity' services that connect the larger stations. 
As we need to serve all stations in the network in each period, we need to keep operating the `Sprinter' lines that stop at all stations. 
Not many changes are made to the frequencies in scenarios where more changes are allowed.
This indicates that the chosen frequency adjustments are maintained in scenarios where more adjustments are allowed.
However, the frequencies of the lines Rtd-Ledn, Rtd-Gvc and Gd-Apn are changed in most scenarios. 
In the scenarios with more changes, the frequency of Rtd-Gvc remains 6 throughout the day, while the frequency of Rtd-Ledn is reduced to 2 in the off-peak scenario. 
These adjustments keep the fast connection between the large stations between The Hague and Rotterdam, while still providing service to the smaller stations on this line. 
The line between Gd and Apn usually takes frequency 2 in the peak periods and 1 in the off-peak, to provide better service in the peak periods. 

When considering the changes made to the stopping pattern, we see a lot of variation in the colours.
This indicates that once certain changes to the stopping pattern are made, they are usually repeated in the scenarios where more adjustments are allowed. 
The adjustments to the stopping pattern are mainly driven by large passenger flows during the peak hour.
For example, the extra stop at Den Haag Laan van NOI (Laa) on the Gvc-Ledn line is often skipped in the peak hours, as there are many passengers who want to travel between Ledn and Gvc during the peak periods. 
By skipping this stop, the in-vehicle time of these passengers is reduced.
However, the passengers travelling between Laa and Ledn are negatively affected by this adjustment, as their in-vehicle time is increased by the time of three additional stops. 
The increased GJT of these passengers turns out to be lower than the reduction in GJT for the passengers travelling between Gvc and Ledn, which is why the adjustment is made.
However, during the off-peak we see fewer adjustments to the stopping pattern, since the small improvement in GJT for some OD pairs does not outweigh the inconvenience for other OD pairs.

\subsubsection{Line plan adjustments with asymmetric lines}
\begin{figure}
    \centering
    \includegraphics[width=0.88\textwidth]{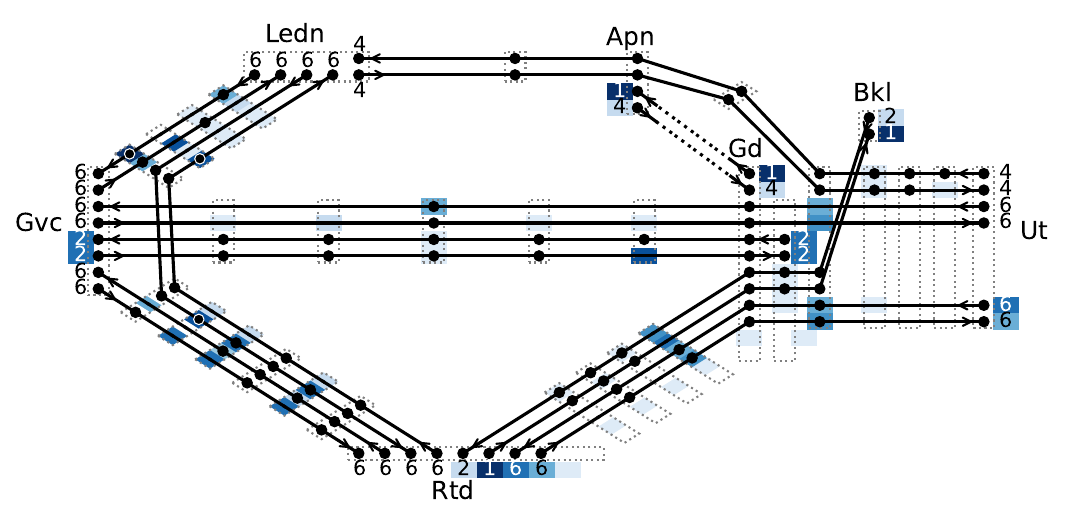}
    \includegraphics[width=0.88\textwidth]{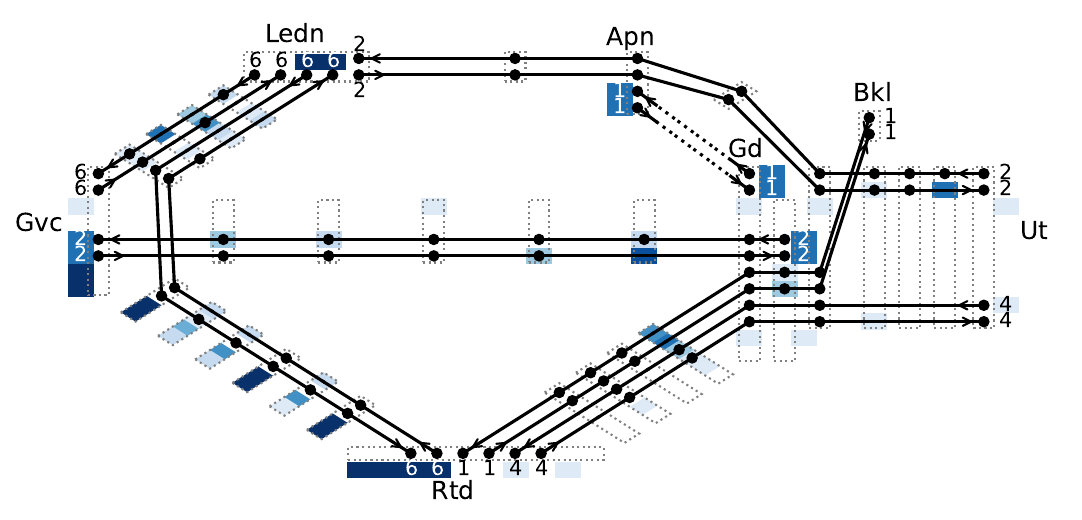}
    \includegraphics[width=0.88\textwidth]{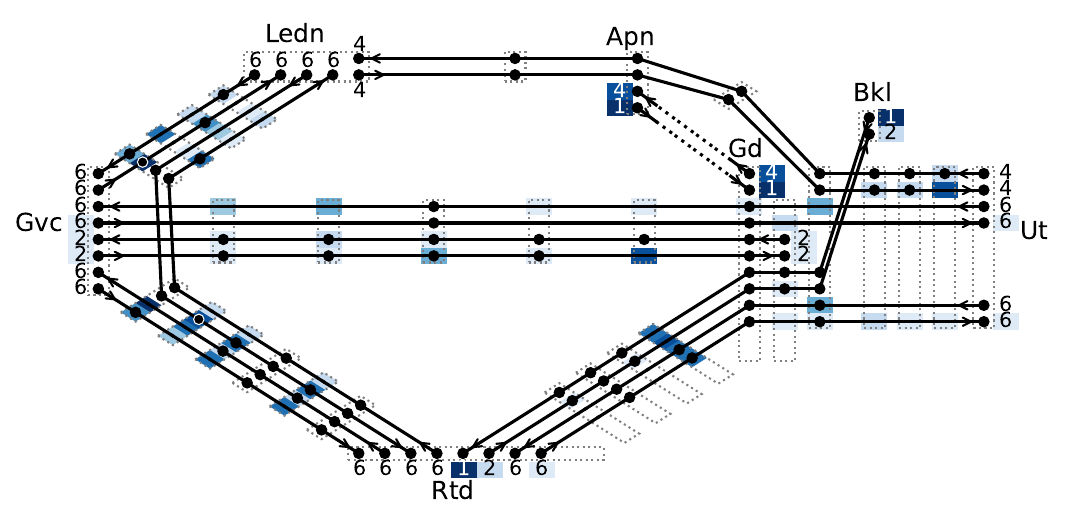}
    \caption{Line plan in \mop, \mo, and \ap\ (from top to bottom) with asymmetric lines and 20 changes. Darker colours denote that the stopping pattern or frequencies were changed more often in scenarios with more changes. The lines between Apn and Gd stop at all intermediate stations in all scenarios.}
    \label{fig:heatmap_asym}
\end{figure}

Similarly to Figure \ref{fig:heatmap_sym}, Figure \ref{fig:heatmap_asym} displays the reference line plan with asymmetric lines and 20 changes.
In this figure, the denoted lines are operated in one direction and the arrows on the line denote this direction.

When we look at the frequencies of the reference line plan, we see that on most routes the frequencies are symmetric and the same as in the reference line plan that uses symmetric lines. 
However, the lines between Rtd and Breukelen (Bkl) and between Apn and Gd do not have the same frequencies in the peak periods, with higher frequencies in the peak direction.
The frequencies on these routes become symmetric during the off-peak. 
Interestingly, these asymmetric frequencies are changed very often and in most scenarios with asymmetric lines the frequencies are the same in both directions. 
For example, in almost all scenarios the lines between Rtd and Bkl have frequency 2 in both directions during the peaks and 1 during the off-peak.
The lines between Apn and Gd have more variation in the frequency, but it is always symmetric in scenarios with more than 20 adjustments. 
For example, in 5 of the 8 additional scenarios both lines take frequency 4 in the morning peak and frequency 2 in the off-peak and afternoon peak.
Another difference compared to frequencies in the symmetric line plans we see on the lines between Gvc and Gouda Goverwelle (Gdgo).
In most scenarios the frequencies on this part of the network are 4 in \mop\ and \mo, and 2 in \ap, compared to a consistent frequency of 2 in the symmetric variants. 
This higher frequency in \mop\ and \mo\ is possible because the two lines between Rtd and Ut take frequency 4 instead of 6 in \mop, and the frequency of the lines between Ledn and Rtd is reduced to 0 instead of 2 in \mo.

Although most frequencies are symmetric, we do see asymmetries in the stopping pattern at some parts of the network, most notably between Gvc and Ledn and between Gvc and Rtd. 
In the reference line plan, the stations Den Haag Mariahoeve (Gvm), Voorschoten (Vst), and De Vink (Dvnk) between Gvc and Ledn are all served by only one line (and hence in only one direction).
Between Gvc and Rtd, most stations are served in both directions.
However, the Intercity services that only serve the large stations and the Sprinter services that serve most or all stations are distributed over different routes. 
In the reference line plan, the lines Gvc-Rtd and Rtd-Ledn serve the large stations, while the lines Rtd-Gvc and Ledn-Rtd stop at most stations between Den Haag HS (Gv) and Rtd.
During the off-peak, both lines between Gvc and Rtd have frequency 0, which results in the line Ledn-Rtd providing the Sprinter service and Rtd-Ledn the Intercity service between Gv and Rtd. 
In several other scenarios, the frequency of both lines between Ledn and Rtd are reduced to 0 in the off-peak, while the frequency of the lines between Gvc and Rtd are kept at 6. 
In those scenarios, the lines between Ledn and Rtd take an Intercity pattern during the peak periods, while the lines between Gvc and Rtd have a Sprinter pattern in both directions. 
During the off-peak, the line Rtd-Gvc keeps this Sprinter pattern, while Gvc-Rtd only stops at the large stations. 
This setup ensures that all stations receive service while facilitating faster connections in one direction.

When comparing the changes made in the scenarios with symmetric lines with the changes made in the scenarios with asymmetric lines there are some similarities.
For example, on most routes the frequencies are reduced by one step in the off-peak, and increased again during the afternoon peak. 
Furthermore, some lines that mostly provide Intercity services (e.g., between Gvc and Ut) are completely removed during the off-peak.
However, by strategically serving some smaller stations with low demand in only one direction, the scenarios with asymmetric lines can achieve a lower GJT. 
By skipping these stations in one direction, faster connections can be provided to OD pairs with larger passenger demand, while keeping some level of service at the stations with lower demand.
Furthermore, as both Sprinter and Intercity services can be combined on lines between Gvc and Rtd by giving each direction a different stopping pattern, the lines between Rtd and Ledn can be completely removed during the off-peak. 
This leaves budget to increase the frequency on the lines between Gvc and Gdgo to 4 instead of 2, which benefits all passengers using this line. 

\section{Conclusion} \label{sec:conclusion}
The quality of the railway service offered to passengers is for a large part determined by the line plan, which influences factors such as direct trips, transfers, and in-vehicle times.
The demand for railway services fluctuates throughout the day, with higher volumes during peak periods and varying demand structures.
As the merits of a line plan depend on the demand it aims to serve, the quality offered by the line plan also varies throughout the day.
Therefore, this paper proposes a model that determines a train line plan with multiple periods that minimises the total GJT. 
Although in recent years there has been more interest in multi-period line planning, most models are proposed for and tested on corridors instead of networks and only incorporate a limited number of service adjustments and symmetric lines.
To address these gaps in the literature, the mixed-integer linear programming model proposed in this paper incorporates the possibility for passengers to transfer, setting stopping patterns and frequencies, and having asymmetric lines. 
The proposed model was tested on a case study based on part of the Dutch railway network. 

Multiple line plans were generated for the case study by varying the allowed number of adjustments to the line plan throughout the day.
Comparing these plans to the reference line plan with symmetric lines and 20 adjustments reveals that allowing for more adjustments can significantly lower the total GJT.
With symmetric lines, additional adjustments can reduce the GJT by up to 1.94\%.
Using asymmetric lines results in even greater reductions, lowering the GJT by 4.26\% compared to the reference scenario.
In general, we see that line plans with asymmetric lines have a lower total GJT than line plans with symmetric lines under the same costs and number of adjustments.
Furthermore, we see that the additional GJT savings decrease when more adjustments are made.
Early adjustments mainly target frequency changes, which impact all passengers of the affected line by altering the waiting and transfer times.
When the frequency of a line is increased, all passengers using that line benefit from it. 
Once frequency adjustments are implemented, the focus shifts to modifying stopping patterns. 
These changes offer less GJT reduction, since they benefit some passengers while inconveniencing others.
In the case with symmetric lines, the changes to the stopping pattern are mainly skipping stops at small to medium large stations that are previously served by two lines.
The majority of the passengers that are using the line do not want to (dis)embark at this station and hence will have a shorter travel time.
In the scenarios with asymmetric lines, stopping pattern adjustments allow for skipping certain stops in only one direction.
This makes it possible to provide a fast service to passengers on large OD pairs, while still providing service to all stations by stopping in one direction. 
While these stopping strategies can reduce the total GJT, only serving a station in one direction has a big impact on passengers who want to travel in the opposite direction from that station.
These passengers will first need to travel farther away from their destination to be able to board a train in the right direction.
How passengers feel about this is not taken into account in this study, which might make the asymmetric line plans look better than they really are. 
Practitioners should weigh both these costs and benefits, to determine if asymmetry in the line plan is beneficial for their network. 

In this paper, we assume that the periods are fixed and part of the input. 
In the case study, demand for the morning hyper peak, midday off-peak and afternoon hyper peak are used, which are all very different. 
One question that would be interesting to tackle in a future study is how many periods can and should be provided to the model. 
When you provide more periods, you can potentially follow the demand better with the line plan.
However, transitioning between periods cannot be done instantly, as trains are driving a certain line somewhere on the network. 
Hence, another direction for future research is how the transition between two line plans can be done efficiently in the corresponding timetable planning.
Lastly, the scalability of the model is also an important aspect for determining the number of periods to provide to the model.
Increasing the network size or the number of periods significantly increases the size of the model.
Therefore, specialised heuristics to solve the resulting model could also be an interesting research direction. 

\section*{Acknowledgements}
\noindent This research is supported by the Netherlands Railways (NS), the largest passenger railway undertaking in the Netherlands.
We also thank Dennis Huisman of NS for his valuable comments and suggestions, which helped to improve this work.

\section*{Disclosure of interest}
The authors report there are no competing interests to declare.

\bibliographystyle{elsarticle-harv}
\bibliography{references}

\newpage
\appendix
\section{Stations in case study} \label{ap:caseStudy}
\begin{table}[ht]
    \centering
    \caption{Overview of station abbreviations used in the case study, E denotes that the station is a terminal station, T denotes that the station is a transfer station.}
    \label{tab:stationAbrreviations}
    \begin{tabular}{llll}
        \hline \Tstrut
        Abbr.          & Station name             & Abbr.           & Station name    \Bstrut\\ \hline \Tstrut
        Apn            & Alphen a/d Rijn (E,T)     & Llzm           & Lansingerland-Zoetermeer \\
        Bdg            & Bodegraven                & Nwk            & Nieuwerkerk a/d IJssel   \\
        Bkl            & Breukelen (E,T)           & Rsw            & Rijswijk                 \\
        Bsk            & Boskoop                   & Rta            & Rotterdam Alexander      \\
        Bsks           & Boskoop Snijdelwijk       & Rtd            & Rotterdam Centraal (E,T) \\
        Cps            & Capelle Schollevaar       & Rtn            & Rotterdam Noord          \\
        Dt             & Delft (T)                 & Sdm            & Schiedam Centrum (T)     \\
        Dtcp           & Delft Campus              & Ut             & Utrecht Centraal (E,T)   \\
        Dvnk           & De Vink                   & Utlr           & Utrecht Leidsche Rijn    \\
        Gd             & Gouda (E,T)               & Utt            & Utrecht Terwijde         \\
        Gdgo           & Gouda Goverwelle (E)      & Vb             & Voorburg                 \\
        Gv             & Den Haag HS (T)           & Vst            & Voorschoten              \\
        Gvc            & Den Haag Centraal (E,T)   & Vtn            & Vleuten                  \\
        Gvm            & Den Haag Mariahoeve       & Wad            & Waddinxveen              \\
        Gvmw           & Den Haag Moerwijk         & Wadn           & Waddinxveen Noord        \\
        Gvy            & Den Haag Ypenburg         & Wadt           & Waddinxveen Triangel     \\
        Laa            & Den Haag Laan van NOI (T) & Wd             & Woerden (T)              \\
        Ldl            & Leiden Lammenschans       & Ztm            & Zoetermeer               \\
        Ledn           & Leiden Centraal (E,T)     & Ztmo           & Zoetermeer Oost          \Bstrut\\ \hline
    \end{tabular}
\end{table}

\end{document}